\documentclass[smallextended]{svjour3}
\usepackage{color}
\usepackage{amssymb}
\usepackage{comment}
\usepackage{graphicx}

\newtheorem{assumption}{Assumption}

\newcommand{\bbH}{\mathbb{H}}

\newcommand{\eqref}[1]{$(\ref{#1})$}
\newcommand{\bproof}{{\it Proof}\,\,}
\newcommand{\eproof}{\hfill$\Box$}

\def\IV{\mathbb{V}}

\def\IP{\mathbb{P}}
\def\IR{\mathbb{R}}

\def\IE{\mathbb{E}}
\def\tV{\mathsf{V}}

\newcommand{\be}{\begin{equation}}
\newcommand{\ee}{\end{equation}}
\newcommand{\beq}{\begin{equation}}
\newcommand{\eeq}{\end{equation}}
\newcommand{\beqas}{\begin{eqnarray*}}
\newcommand{\eeqas}{\end{eqnarray*}}

\newcommand{\bbR}{\mathbb R}

\newcommand{\cG}{{\mathcal G}}
\newcommand{\cF}{{\mathcal F}}

\definecolor{darkred}{rgb}{.7,0,0}

\definecolor{darkgreen}{rgb}{0,0.7,0}

\definecolor{darkblue}{rgb}{0,0,0.7}

\newcommand{\bea}{\begin{eqnarray}}
\newcommand{\eea}{\end{eqnarray}}
\newcommand{\beas}{\begin{eqnarray*}}
\newcommand{\eeas}{\end{eqnarray*}}

\def\IP{\mathbb{P}}
\def\IR{\mathbb{R}}

\def\IH{\mathbb{H}}
\def\IX{\mathbb{X}}

\def\cH{{\mathcal H}}

\def\cF{{\mathcal F}}
\def\cG{{\mathcal G}}

\def\cP{{\mathcal P}}
\def\cQ{{\mathcal Q}}

\def\cV{{\mathbb X}}
\def\bbX{{\mathbb X}}
\newcommand{\sQ}{\mathsf Q}

\begin{document}
\title{Determining White Noise Forcing From Eulerian
Observations in the Navier-Stokes Equation}

\titlerunning{Determining White Noise Forcing of the Navier Stocks Equation} 
\author{Viet Ha Hoang \and Kody J. H. Law \and Andrew~M.~Stuart}

\institute{
Viet Ha Hoang \at Division of Mathematical Sciences, School of Physical and Mathematical Sciences, Nanyang Technological
University, Singapore, 637371.
\email{vhhoang@ntu.edu.sg}
\and
Kody J. H. Law \at Mathematics Institute, Warwick University, Coventry CV4 7AL, UK.
\email {k.j.h.law@warwick.ac.uk}
\and Andrew M. Stuart \at Mathematics Institute, Warwick University, Coventry CV4 7AL, UK.
\email{a.m.stuart@warwick.ac.uk}
}



\maketitle
\begin{abstract}

The Bayesian approach to inverse problems is of paramount
importance in quantifying uncertainty about the input to, and the state 
of, a system of interest given noisy observations. Herein we 
consider the forward problem of the forced 2D Navier-Stokes equation.
The inverse problem is to make inference concerning the forcing, and possibly the
initial condition, given noisy observations
of the velocity field.  We place a prior on the forcing which is in the
form of a spatially-correlated and temporally-white Gaussian process,
and formulate the inverse problem for the posterior distribution.
Given appropriate spatial regularity
conditions, we show that the solution
is a continuous function of the forcing.  Hence, for 
appropriately chosen spatial regularity in the prior,  
the posterior distribution on the forcing is absolutely continuous 
with respect to the prior and is hence well-defined.  Furthermore,
it may then be shown that
the posterior distribution is a continuous function of the data.
We complement these theoretical results with numerical simulations
showing the feasibility of computing the posterior distribution, and
illustrating its properties.

\end{abstract}

\section{Introduction} \label{sec:I}

The Bayesian approach to inverse problems has grown in popularity
significantly over the last decade, driven by algorithmic innovation
and steadily increasing computer power \cite{kaipio2004statistical}.
Recently there have been systematic developments of the theory
of Bayesian inversion on function space \cite{La02,La07,CDRS,stuart2010inverse,lasanen2012non,lasanen2012nonb} 
and this has led to new sampling algorithms which perform well
under mesh-refinement \cite{CRSW12,vollmer2013,law2012proposals}. 
In this paper we add to this growing interest in the Bayesian formulation of
inversion, in the context of a specific PDE inverse problem, 
motivated by geophysical applications such as 
data assimilation in the atmosphere and ocean
sciences, and demonstrate that
fully Bayesian probing of the posterior distribution is feasible.

The primary goal of this paper is to demonstrate 
that the Bayesian formulation
of inversion for the forced Navier-Stokes equation, introduced in
\cite{CDRS}, can be extended to the case of white noise forcing.
The paper \cite{CDRS} assumed an Ornstein-Uhlenbeck structure in 
time for the forcing, and hence did not include the white noise case. 
It is technically demanding to extend to the case of white noise forcing, 
but it is also of practical interest. This practical importance stems 
from the fact that the Bayesian formulation of problems with white noise forcing
corresponds to a statistical version of the continuous time weak
constraint 4DVAR methodology \cite{zupanski1997general}.
The 4DVAR approach to data assimilation currently gives the most accurate
global short term weather forecasts available \cite{lorenc2003potential}
and this is arguably the case because, unlike ensemble filters which
form the major competitor, 4DVAR has a rigorous statistical interpretation
as a maximum a posteriori (or MAP) estimator -- the point
which maximizes the posterior probability.
It is therefore of interest to seek to embed our
understanding of such methods in a broader Bayesian context.

To explain the connection between our work and the 4DVAR methodology
and, in particular, to explain the terminology used in the
data assimilation community,
it is instructive to consider the finite dimensional differential
equation
$$\frac{du}{dt}=f(u)+\xi, \quad u(0)=u_0$$
on $\bbR^n$.  Assume that we are given noisy observations $\{y_j\}$
of the solution $u_j=u(t_j)$ at times $t_j=jh$ so that
$$y_j=u_j+\eta_j, \quad j=1,\cdots, J,$$
for some noises $\eta_j.$ An important inverse problem is to
determine the initial condition $u_0$ and forcing $\xi$
which best fits the data. If we view the solution $u_j$ as
a function of the initial condition and forcing,
then a natural regularized least squares problem is to determine 
$u_0$ and $\xi$ to minimize
$$I(u_0,\xi)=\sum_{j=1}^J \Bigl|\Gamma^{-\frac12}\bigl(y_j-u_j(u_0,\xi)
\bigr)\Bigr|^2+|\Sigma^{-\frac12}u_0|^2+\|\sQ^{-\frac12}\xi\|^2$$
where $|\cdot|, \|\cdot\|$ denote the $\bbR^n-$Euclidean and $L^2\bigl([0,T];\bbR^n\bigr)$
norms respectively, $\Gamma, \Sigma$ denote covariance matrices
and $\sQ$ a covariance operator. This is a continuous time analogue
of the 4DVAR or variational methodology, as described in 
the book of Bennett \cite{bennett2002inverse}. In numerical
weather prediction the method is known as weak constraint 4DVAR,
and as 4DVAR if $\xi$ is set to zero (so that the ODE
$$\frac{du}{dt}=f(u), \quad u(0)=u_0$$
is satisfied as a hard constraint), the term $\|\sQ^{-\frac12}\xi\|^2$
dropped, and the minimization is over $u_0$ only.
As explained in \cite{kaipio2004statistical}, these minimization
problems can be viewed as probability maximizers for the posterior
distribution of a Bayesian formulation of the inverse problem --
the so-called MAP estimators. In this interpretation the
prior measures on $u_0$ and $\xi$ are centred Gaussians with
covariances $\Sigma$ and $\sQ$ respectively. Making this connection opens up the
possibility of performing rigorous statistical inversion, and thereby
estimating uncertainty in predictions made.

The ODEs arising in atmosphere and ocean science applications are of very high
dimension, arising from discretizations of PDEs. It is therefore
conceptually important to carry through the program in the previous
paragraph, and in particular Bayesian formulation of the inversion,
for PDEs; 
the paper \cite{dashti2013map} explains how to define
MAP estimators for measures on Hilbert spaces and the connection to
variational problems.
The Navier-Stokes equation in 2D provides a useful canonical example 
of a PDE of direct relevance to the atmosphere and ocean sciences.
{When the prior covariance operator $\sQ$ is chosen to be that associated
to an Ornstein-Uhleneck operator in time, the Bayesian formulation
for the 2D Navier-Stokes equation has been carried out in \cite{CDRS}.
Our goal in this paper is to extend to the more technically demanding
case where $\sQ$ is the covariance operator associated with
a white noise in time, with spatial
correlations $Q$. 
We will thus use the prior model $\xi dt={dW}$ where $W$ is a $Q-$Wiener process in
an appropriate Hilbert space, and consider
inference with respect to $W$ and $u_0.$
In the finite dimensional setting the differences between the case
of coloured and white noise forcing, with respect to the inverse
problem, are much less substantial and the interested reader
may consult \cite{HSV11} for details.
}

The key tools required in applying the function space Bayesian approach
in \cite{stuart2010inverse} are the proof of continuity of the forward map
from the function space of the unknowns to the data space, together
with estimates of the dependence of the forward map upon its point
of application, sufficient to show certain integrability properties with
respect to the prior. This program is carried out for the 2D Navier-Stokes
equation with Ornstein-Uhlenbeck priors on the forcing in the paper
\cite{CDRS}. However to use priors which are white in time adds
further complications since it is necessary to study the stochastically
forced 2D Navier-Stokes equation and to establish continuity of the
solution with respect to small changes in the Brownian motion $W$ which defines
the stochastic forcing. We do this by employing the solution concept
introduced by Flandoli in \cite{Flandoli94}, and using probabilistic
estimates on the solution derived by Mattingly in \cite{Mattingly99}.
In section \ref{sec:PS} we describe the relevant theory of the forward
problem, employing the setting of Flandoli. In section \ref{sec:B} we
build on this theory, using the estimates of Mattingly to verify the
conditions in \cite{stuart2010inverse}, resulting in a well-posed Bayesian
inverse problem for which the posterior is Lipschitz in the data
with respect to Hellinger metric. Section \ref{sec:C} extends this
to include making inference about the initial condition as well
as the forcing. Finally, in section \ref{sec:D}, we present numerical
results which demonstrate feasibility of sampling from the
posterior on white noise forces, {and demonstrate the properties
of the posterior distribution.}

\section{Forward Problem}
\label{sec:PS}

In this section we study the forward problem of the Navier-Stokes
equation driven by white noise. Subsection
\ref{ssec:over} describes the forward problem, the
Navier-Stokes equation, and rewrites it as an 
ordinary differential equation in a Hilbert space.
In subsection \ref{ssec:spaces} we define the
functional setting used throughout the paper.
Subsection \ref{ssec:solc} highlights the solution
concept that we use, leading 
in subsection \ref{ssec:cont} to proof of 
the key fact that the solution of the Navier-Stokes
equation is continuous as a function of the
rough driving of interest and the initial condition.
{All our theoretical results in this paper are derived
in the case of Dirichlet (no flow) boundary
conditions.  They may be extended 
to the problem on the periodic torus $\mathbb{T}^d$, 
but we present the more complex Dirichlet case only for brevity.}

\subsection{Overview}
\label{ssec:over}
Let $D\in\IR^2$ be a bounded domain with smooth boundary. 
We consider in $D$ the Navier-Stokes equation 
\begin{eqnarray}
\partial_t u-\nu\Delta u+u\cdot\nabla u&=&f-\nabla p,\quad (x,t) \in D \times (0,\infty)\label{eq:NSE}\\
\nabla\cdot u&=&0,\quad (x,t) \in D \times (0,\infty)\nonumber\\
u&=&0,\quad (x,t)\in\partial D \times (0,\infty),\nonumber\\
u&=&u_0,\quad (x,t) \in D \times \{0\}.\nonumber 
\end{eqnarray}
We assume that the initial condition $u_0$ and
the forcing $f(\cdot,t)$ are divergence-free. We will in particular work with equation \eqref{eq:u} below, obtained by projecting \eqref{eq:NSE} into the space of divergence free functions -- {the Leray projector \cite{temam1984navier}.}
We denote by $\tV$ the space of all 
divergence-free smooth functions from $D$ to $\IR^2$ with
compact support, by $\IH$ the closure of $\tV$ in $(L^2(D))^2$, 
and by $\IH^1$ the closure of $\tV$ in $(H^1(D))^2$. 
Let $\IH^2=(H^2(D))^2\bigcap\IH^1$.  The initial condition $u_0$ is assumed to be in $\IH$. We define the linear 
Stokes' operator $A:\IH^2\to \IH$ by $A u=-\Delta u$  noting that the assumption
of compact support means that Dirichlet boundary condition are
imposed on the Stokes' operator $A$. 
Since $A$ is selfadjoint, $A$ possesses eigenvalues 
$0<\lambda_1\le\lambda_2\le\ldots$ with the corresponding 
eigenvectors $e_1, e_2, \ldots\in \IH^2$.

We denote by $\langle\cdot,\cdot\rangle$ the inner product 
in $\IH$,  extended to the dual pairing on $\IH^{-1} \times \IH^1$.
We then define the bilinear form $B: \IH^1\times\IH^1\to \IH^{-1}$
\[
\langle B(u,v),z\rangle=\int_Dz(x)\cdot(u(x)\cdot\nabla)v(x)dx
\]
which must hold for all $z \in \IH^1.$
From the incompressibility condition we have, {for
all $z \in \IH^1$},
\be
\langle B(u,v),z\rangle=-\langle B(u,z),v\rangle.
\label{eq:Buvz}
\ee
By projecting problem \eqref{eq:NSE} into $\IH$ 
we may write it as an ordinary differential equation
in the form 
\be
du(t)=-\nu A udt-B(u,u)dt+dW(t),\ \ u(0)=u_0\in\IH,
\label{eq:u}
\ee 
where $dW(t)$ is the projection of the forcing $f(x,t)dt$ 
into $\IH$. We will define the solution of this equation
pathwise, for suitable $W$, {not necessarily differentiable
in time}. 

\subsection{Function Spaces}
\label{ssec:spaces}

For any $s\ge 0$ we define
$\IH^s\subset\IH$ to be the Hilbert space of functions $u=\sum_{k=1}^\infty u_ke_k\in \IH$ such that
\[
\sum_{k=1}^\infty\lambda_k^{s}u_k^2<\infty;
\]
we note that the $\IH^j$ for $j \in \{0,1,2\}$ coincide
with the preceding definitions of these spaces. 
The space $\IH^s$ is endowed with the inner product
\[
\langle u,v\rangle_{\IH^s}=\sum_{k=1}^\infty\lambda_k^{s}u_kv_k,
\]
{for $u=\sum_{k=1}^\infty u_ke_k$, $v=\sum_{k=1}^\infty v_ke_k$
in $\IH$}.
We denote by $\IV$ the particular choice $s=\frac12+\epsilon$,
namely $\IH^{\frac12+\epsilon}$, for given $\epsilon>0$. 
In what follows we will be particularly interested in
continuity of the mapping from the forcing  
$W$ into linear functionals of the solution of
\eqref{eq:u}. To this end it is helpful to
define the Banach space $\cV:=C([0,T];\IV)$ with the norm
\[
\|W\|_{\cV}=\sup_{t\in (0,T)}\|W(t)\|_{\IV}.
\]

\subsection{Solution Concept}
\label{ssec:solc}

In what follows we define a solution concept
for equation \eqref{eq:u} for each forcing function
$W$ which is continuous, but not necessarily differentiable,
in time. We always assume that $W(0)=0.$ Following 
Flandoli \cite{Flandoli94}, for each $W\in \cV$ we define the weak solution 
$u(\cdot;W)\in C([0,T];\IH)\bigcap L^2([0,T];\IH^{1/2})$  
of \eqref{eq:u} as a function that satisfies
\be
\langle u(t),\phi\rangle+{\nu}\int_0^t\langle u(s),A\phi\rangle ds-\int_0^t\langle {B\bigl(u(s),\phi\bigr),u(s)}\rangle dx=
\langle u_0,\phi\rangle+\langle W(t),\phi\rangle,
\label{eq:weaksln}
\ee
for all $\phi\in \IH^2$ and all $t\in (0,T)$; note the
integration by parts on the Stokes' operator and the
use of \eqref{eq:Buvz} to derive this identity from \eqref{eq:u}. Note further that if $u$ and $W$ are sufficiently smooth, \eqref{eq:weaksln} is equivalent to \eqref{eq:u}.

{To employ this solution concept we first introduce the concept of a solution 
of the linear equation}
\be
dz(t)=-\nu A zdt+dW(t),\ \ z(0)=0\in\IH
\label{eq:z}
\ee 
where $W$ is a deterministic continuous function obtaining values in $\IX$ but not necessarily differentiable.
We define a weak solution of this equation as a function $z\in C([0,T];\IH)$ such that
\be
\langle z(t),\phi\rangle+{\nu} \int_0^t\langle z(s),A\phi\rangle ds=\langle W(t),\phi\rangle
\label{eq:weakz}
\ee
for all $\phi\in\IH^2$.

Then for this function $z(t)$ we consider the solution $v$ of
the equation
\be
dv(t)=-\nu A vdt-B(z+v,z+v)dt,\ \ v(0)=u_0\in\IH.
\label{eq:v}
\ee 
As we will show below, $z(t)$ possesses sufficiently regularity so \eqref{eq:v} possesses a weak solution $v$. We then 
deduce that $u=z+v$ is a weak solution of \eqref{eq:u} in the sense of 
\eqref{eq:weaksln}. When we wish to emphasize the dependence of $u$ on $W$
(and similarly for $z$ and $v$) we write $u(t;W).$

We will now show that the function $z$ defined by
\begin{eqnarray}
z(t)&=&\int_0^t e^{-\nu A(t-s)}dW(s)\nonumber\\
&=&W(t)-\int_0^ t\nu A e^{-\nu A(t-s)}W(s)ds.
\label{eq:zt}
\end{eqnarray}
satisfies the weak formula \eqref{eq:weakz}. Let $w_k=\langle W, e_k \rangle$, that is
\be
\label{eq:W}
W(t):=\sum_{k=1}^\infty {w_k}(t)e_k\in \cV.
\ee
We then deduce from \eqref{eq:zt} that
\be
z(t;W)=W(t)-\sum_{k=1}^\infty \Bigl(\int_0^t w_k(s){\nu}\lambda_ke^{(t-s)(-{\nu}\lambda_k)}ds\Bigr)e_k.
\label{eq:zz}
\ee
We have the following regularity property for $z$:
\begin{lemma}
\label{lem:z2}
For each $W \in\cV $, the
function $z=z(\cdot;W)\in C([0,T];\IH^{1/2})$.
\end{lemma}

\bproof
We first show that for each $t$, $z(t;W)$ as defined in \eqref{eq:zz} belongs to $\IH^{1/2}$. Fixing an integer $M>0$, using inequality $a^{1-\epsilon/2}e^{-a}<c$ for all $a>0$ for an appropriate constant $c$, we have
\beqas
\sum_{k=1}^M\lambda_k^{1/2}\Bigg(\int_0^t\nu\lambda_ke^{(t-s)(-\nu\lambda_k)}w_k(s)ds\Bigg)^2
\le\sum_{k=1}^M\lambda_k^{1/2}\Bigg(\int_0^t{c\over (t-s)^{1-\epsilon/2}}\lambda_k^{\epsilon/2}|w_k(s)|dx\Bigg)^2.
\eeqas
Therefore,
\beqas
\Bigg\|\sum_{k=1}^M\int_0^t\nu\lambda_ke^{(t-s)(-\nu\lambda_k)}w_k(s)e_kds\Bigg\|_{\IH^{1/2}}
&\le& \Bigg\|\sum_{k=1}^M\int_0^t{c\over(t-s)^{1-\epsilon/2}}\lambda_k^{\epsilon/2}|w_k(s)|e_kds\Bigg\|_{\IH^{1/2}}\\
&\le& \int_0^t{c\over(t-s)^{1-\epsilon/2}}\Bigg\|\sum_{k=1}^M\lambda_k^{\epsilon/2}|w_k(s)|e_k\Bigg\|_{\IH^{1/2}}ds\\
&\le&\max_{s\in(0,T)}\|W(s)\|_{\IH^{1/2+\epsilon}}\int_0^t{c\over (t-s)^{1-\epsilon/2}}ds,
\eeqas
which is uniformly bounded for all $M$. Therefore, 
\[
\sum_{k=1}^{\infty}\Bigg(\int_0^tw_k(s)\nu\lambda_ke^{(t-s)(-\nu\lambda_k)}ds\Bigg)e_k\in \IH^{1/2},
\]
It follows from \eqref{eq:zz} that, since $W \in \cV$,
for each $t$, $z(t;W)\in \IH^{1/2}$ as required. Furthermore, for all $t\in (0,T)$
\be
\|z(t;W)\|_{\IH^{1/2}}\le c\|W\|_{\cV}.
\label{eq:boundH12ofz}
\ee
Now we turn to the continuity in time.
Arguing similarly, we have that
\beqas
&&\Bigg\|\sum_{k=M}^\infty\Bigg(\int_0^tw_k(s)\nu\lambda_ke^{(t-s)(-\nu\lambda_k)}ds\Bigg)e_k\Bigg\|_{\IH^{1/2}}
\le\int_0^t{c\over (t-s)^{1-\epsilon/2}}\Bigg \|\sum_{k=M}^\infty w_k(s)e_k\Bigg\|_{\IH^{1/2+\epsilon}}ds\\
&&\qquad\qquad\le \Bigg(\int_0^t{c\over (t-s)^{(1-\epsilon/2)^p}}ds\Bigg)^{1/p}\Bigg(\int_0^t\Bigg \|\sum_{k=M}^\infty w_k(s)e_k\Bigg\|^q_{\IH^{1/2+\epsilon}}ds \Bigg)^{1/q},
\eeqas
for all $p,q>0$ such that $1/p+1/q=1$. From 
the Lebesgue dominated convergence theorem,
\[
\lim_{M\to\infty}\int_0^{t}\Bigg\|\sum_{k=M}^\infty w_k(s)e_k\Bigg\|^q_{\IH^{1/2+\epsilon}}ds=0;
\]
and when $p$ sufficiently close to 1, 
\[
\int_0^t{c\over (t-s)^{(1-\epsilon/2)p}}ds
\]
is finite. We then deduce that
\[
\lim_{M\to\infty}\Bigg\|\sum_{k=M}^\infty\Bigg(\int_0^tw_k(s)\nu\lambda_ke^{(t-s)(-\nu\lambda_k)}ds\Bigg)e_k\Bigg\|_\IH^{1/2}=0,
\]
uniformly for all $t$. 

Fixing $t\in (0,T)$ we show that
\[
\lim_{t'\to t}\|z(t;W)-z({t';W})\|_{\IH^{1/2}}=0.
\]
We have
\beqas
&&\|z(t;W)-z(t';W)\|_{\IH^{1/2}}\le \|W(t)-W(t')\|_{\IH^{1/2}}+\\
&&\qquad\qquad\Bigg\|\sum_{k=1}^{M-1}\Bigg(\int_0^tw_k(s)\nu\lambda_ke^{(t-s)(-\nu\lambda_k)}ds-\int_0^{t'}w_k(s)\nu\lambda_ke^{(t'-s)(-\nu\lambda_k)}ds\Bigg)e_k\Bigg\|_{\IH^{1/2}}+\\
&&\Bigg\|\sum_{k=M}^\infty\Bigg(\int_0^tw_k(s)\nu\lambda_ke^{(t-s)(-\nu\lambda_k)}ds\Bigg)e_k\Bigg\|_{\IH^{1/2}}+\Bigg\|\sum_{k=M}^\infty\Bigg(\int_0^{t'}w_k(s)\nu\lambda_ke^{(t'-s)(-\nu\lambda_k)}ds\Bigg)e_k\Bigg\|_{\IH^{1/2}}.
\eeqas
For $\delta>0$, when $M$ is sufficiently large, the argument above shows that
\[
\Bigg\|\sum_{k=M}^\infty\Bigg(\int_0^tw_k(s)\nu\lambda_ke^{(t-s)(-\nu\lambda_k)}ds\Bigg)e_k\Bigg\|_{\IH^{1/2}}+\Bigg\|\sum_{k=M}^\infty\Bigg(\int_0^{t'}w_k(s)\nu\lambda_ke^{(t'-s)(-\nu\lambda_k)}ds\Bigg)e_k\Bigg\|_{\IH^{1/2}}<\delta/3. 
\]
Furthermore, when $|t'-t|$ is sufficiently small,
\[
\Bigg\|\sum_{k=1}^{M-1}\Bigg(\int_0^tw_k(s)\nu\lambda_ke^{(t-s)(-\nu\lambda_k)}ds-\int_0^{t'}w_k(s)\nu\lambda_ke^{(t'-s)(-\nu\lambda_k)}ds\Bigg)e_k\Bigg\|_{\IH^{1/2}}<\delta/3.
\]
Finally, since $W \in \cV$, for $|t'-t|$ is sufficiently small we have
\[
\|W(t)-W(t')\|_{\IH^{1/2}} <\delta/3
\]
Thus when $|t'-t|$ is sufficiently small, $\|z(t;W)-z(t';W)\|_{\IH^{1/2}}<\delta$. The conclusion follows.
\eproof 

Having established regularity, we now show that $z$ is indeed
a weak solution of \eqref{eq:z}.

\begin{lemma}
\label{lem:z}
For each $\phi\in\IH^2$, $z(t)=z(t;W)$ satisfies \eqref{eq:weakz}.
\end{lemma}
\bproof\ 
It is sufficient to show this for $\phi=e_k$. We have
\beqas
\int_0^t\langle z(s),Ae_k\rangle ds=\int_0^t\langle W(s),Ae_k\rangle ds-\int_0^t\int_0^s {w_k(\tau)\nu}\lambda_k^2e^{(s-\tau)(-
{\nu}\lambda_k)}d\tau ds\\
=\lambda_k\int_0^t {w_k}(s)ds-{\nu}\lambda_k^2\int_0^t
{w_k}(\tau)\Bigl(\int_\tau^te^{(s-\tau)(-{\nu}\lambda_k)}
ds\Bigr)d\tau\\
=\lambda_k\int_0^t{w_k}(s)ds-\lambda_k\int_0^t{w_k}(\tau)d\tau+\lambda_k\int_0^t{w_k}(\tau)e^{(t-\tau)(-{\nu}\lambda_k)}d\tau\\
=\lambda_k\int_0^t{w_k}(\tau)e^{(t-\tau)(-{\nu}\lambda_k)}d\tau.
\eeqas
On the other hand,
\[
\langle z(t),e_k\rangle =\langle W(t),e_k\rangle -{\nu}
\lambda_k\int_0^t{w_k}(s)e^{(t-s)(-\nu\lambda_k)}ds.
\]
The result then follows.
\eproof

We now turn to the following result, which
concerns $v$ and is established on page 416 of
\cite{Flandoli94}, given 
the properties of $z(\cdot;W)$ established
in the preceding two lemmas.

\begin{lemma}\label{lem:v}
For each $W\in\cV$, problem \eqref{eq:v} has a unique solution 
$v$ in the function space $C(0,T;\IH)\bigcap L^2(0,T;\IH^1)$.
\end{lemma}

{We then have the following existence and uniqueness result
for the Navier-Stokes equation \eqref{eq:u},
more precisely for the weak form \eqref{eq:weaksln}, driven
by rough additive forcing \cite{Flandoli94}:} 

\begin{proposition}
\label{prop:u}
For each $W\in\cV$,
problem \eqref{eq:weaksln} has a unique solution $u\in C(0,T;\IH)\bigcap L^2(0,T;\IH^{1/2})$ such that $u-z\in L^2(0,T;\IH^1)$.
\end{proposition}
\bproof
A solution $u$ for \eqref{eq:weaksln} can be taken as
\be
u(t;W)=z(t;W)+v(t;W).
\label{eq:defu}
\ee
From the regularity properties of $z$ and $v$ in Lemmas \ref{lem:z2} and \ref{lem:v}, we deduce that $u\in C(0,T;\IH)\bigcap L^2(0,T;\IH^{1/2})$.
Assume that $\bar u(t;W)$ is another solution of \eqref{eq:weaksln}. Then $\bar v(t;W)=\bar u(t;W)-z(t;W)$ is a solution in $C(0,T;\IH)\bigcap L^2(0,T;\IH^1)$ of \eqref{eq:v}. However, \eqref{eq:v} has a unique solution in $C(0,T;\IH)\bigcap L^2(0,T;\IH^1)$. Thus $\bar v=v$.
\eproof

\subsection{Continuity of the Forward Map}
\label{ssec:cont}

The purpose of this subsection is to establish
continuity of the forward map from $W$ into the weak solution
$u$ of \eqref{eq:u}, {as defined in \eqref{eq:weaksln},}
at time $t>0.$
In fact we prove continuity of the forward map from
$(u_0,W)$ into $u$ and for this it is useful to  
define the space $\cH=\IH\times\cV$
and denote the solution $u$ by $u(t;u_0,W)$.

\begin{theorem}
\label{t:first}
For each $t>0$, the solution $u(t;\cdot,\cdot)$ of \eqref{eq:u}
is a continuous map from $\cH$ into $\IH$. 
\end{theorem}

\bproof 
First we fix the initial condition and just write $u(t;W)$ for
simplicity.
We consider equation \eqref{eq:u} with driving $W \in \IX$ given by
\eqref{eq:W} and by $W' \in \IX$ defined by
\[
W'(s)=\sum_{k=1}^\infty {w_{k}}'(s)e_{k}\in\cV.
\]
We will prove that, {for $W,W'$ from a bounded set
in $\IX$}, there is $c=C(T)>0$, such that 
\be
\label{eq:zc}
 \sup_{t\in(0,T)}\|z(t;W)-z(t;W')\|_{\IH^{1/2}}\le c \|W-W'\|_{\cV}
\ee
and, {for each $t \in (0,T)$}, 
\be
\label{eq:vc}
\|v(t;W)-v(t;W')\|_{\IH}^2\le c\sup_{s\in(0,T)}\|z(s;W)-z(s;W')\|_{L^4(D)}^2. 
\ee
This suffices to prove the desired result since Sobolev
embedding yields, from \eqref{eq:vc},
\be
\label{eq:vc2}
\|v(t;W)-v(t;W')\|_{\IH}^2\le c\sup_{s\in(0,T)}\|z(s;W)-z(s;W')\|_{\bbH^{\frac12}}^2. 
\ee
Since $u=z+v$ we deduce from \eqref{eq:zc} and \eqref{eq:vc2}
that $u$ as a map from $\cV$ to $\IH$ is continuous.

To prove \eqref{eq:zc} we note that 
\[
\|z(t;W)-z(t;W')\|_{\bbH^{\frac12}} \le \|W(t)-W'(t)\|_{\bbH^{\frac12}}+
\Bigg\| \int_0^t \nu Ae^{-\nu A(t-s)}\bigl(W(s)-W'(s)\bigr)ds\Bigg\|_{\bbH^{\frac12}}
\]
so that
\beqas
&&\sup_{t \in(0,T)}\|z(t;W)-z(t;W')\|_{\bbH^{\frac12}} \le
\|W-W'\|_{\cV}\\
&&\quad\quad\quad\quad\quad\quad\quad\quad\quad\quad\quad\quad\quad\quad\quad+\sup_{t \in (0,T)}\Bigg\| \int_0^t \nu Ae^{-\nu A(t-s)}\bigl(W(s)-W'(s)\bigr)ds\Bigg\|_{\bbH^{\frac12}}.
\eeqas
Thus it suffices to consider the last term on the right hand side. We have
\beqas
&&\Bigg\| \int_0^t  Ae^{-\nu
    A(t-s)}\bigl(W(s)-W'(s)\bigr)ds\Bigg\|_{\bbH^{\frac12}}^2\\
&&\qquad=\Bigg\|\sum_{k=1}^{\infty}\int_0^t\lambda_ke^{(t-s)(-\nu\lambda_k)}(w'_k(s)-w_k(s))e_kds\Bigg\|_{\IH^{1/2}}^2\\
&&\qquad=\sum_{k=1}^\infty\lambda_k^{1/2}\Bigg(\int_0^t\lambda_ke^{(t-s)(-{\nu}\lambda_k)}(w_k'(s)-w_k(s))ds\Bigg)^2\\
&&\qquad\le\sum_{k=1}^\infty\lambda_k^{1/2}\Bigg(\int_0^t\lambda_ke^{(t-s)(-{\nu}\lambda_k)}|w_k'(s)-w_k(s)|ds\Bigg)^2\\
&&\qquad\le\sum_{k=1}^\infty\lambda_k^{1/2}\Bigg(\int_0^t{c\over (t-s)^{1-\epsilon/2}}\lambda_k^{\epsilon/2}|w_k'(s)-w_k(s)|ds\Bigg)^2
\eeqas
where we have used the fact that $a^{1-\epsilon/2}e^{-a}<c$ for all $a>0$ for an appropriate constant $c$. From this, we deduce
that
\beqas
&&\Bigg\| \int_0^t  Ae^{-\nu A(t-s)}\bigl(W(s)-W'(s)\bigr)ds\Bigg\|_{\bbH^{\frac12}}\\
&&\qquad\le\Bigg\|\sum_{k=1}^\infty\int_0^t{c\over (t-s)^{1-\epsilon/2}}\lambda_k^{\epsilon/2}|w_k'(s)-w_k(s)|e_kds\Bigg\|_{\IH^{1/2}}\\
&&\qquad\le\int_0^t{c\over (t-s)^{1-\epsilon/2}}\Bigg\|\sum_{k=1}^\infty\lambda_k^{\epsilon/2}|w_k'(s)-w_k(s)|e_k\Bigg\|_{\IH^{1/2}}ds\\
&&\qquad\le\int_0^t{c\over (t-s)^{1-\epsilon/2}}ds\sup_{s\in(0,T)}\Bigg\|\sum_{k=1}^\infty\lambda_k^{\epsilon/2}|w_k'(s)-w_k(s)|e_k\Bigg\|_{\IH^{1/2}}\\
&&\qquad\le c\sup_{s\in(0,T)}\|W'(s)-W(s)\|_{\IV}.
\eeqas
Therefore \eqref{eq:zc} holds.

We now prove \eqref{eq:vc}.
We will use the following estimate for the
solution $v$ of \eqref{eq:v} which is proved in 
Flandoli \cite{Flandoli94}, page 412, by means
of a Gronwall argument:
\be
\sup_{s\in(0,T)} \|v(s)\|^2_{\IH}+\int_0^T\|v(s)\|^2_{\IH^1}\le C(T,\sup_{s\in(0,T)}\|z(s)\|_{L^4(D)}).
\label{eq:bound}
\ee
We show that the map $C([0,T];L^4(D))\ni z(\cdot;W)\mapsto v(\cdot;W)\in\IH$
is continuous.
For $W$ and $W'$ in $\cV$, {define $v=v(t;W), v'=v(t;W'),
z=z(t;W), z'=z(t;W'), e=v-v'$
and $\delta=z-z'$.} Then we have 
\be
{de\over dt}+\nu Ae+
B\bigl(v+z,v+z\bigr)-B\bigl(v'+z',v'+z'\bigr)=0.
\ee
From this, we have
\beqas
&&{1\over 2}{d\|e\|_{\IH}^2\over dt}+\nu\|e\|_{\IH^1}^2=\\
&&\quad-\bigl\langle B\bigl(v+z,v+z\bigr),e\bigr\rangle\\
&&\quad\quad+\big\langle B\bigl(v'+z',v'+z'\bigr),e\big\rangle.
\eeqas
From \eqref{eq:Buvz} we obtain
\beqas
&&{1\over 2}{d\|e\|_{\IH}^2\over dt}+\nu\|e\|_{\IH^1}^2=\\ 
&&\quad+\bigl\langle B\bigl(v+z,e\bigr),v+z\bigr\rangle\\
&&\quad\quad-\big\langle B\bigl(v'+z',e\bigr),v'+z'\big\rangle\\
&&\quad=\bigl\langle B\bigl(v+z,e\bigr),v+z-v'-z'\bigr\rangle\\
&&\quad\quad-\big\langle B\bigl(v'+z'-v-z,e\bigr),v'+z'\big\rangle\\
&&\quad=\bigl\langle B\bigl(v+z,e\bigr),e+\delta\bigr\rangle\\
&&\quad\quad+\big\langle B\bigl(e+\delta,e\bigr),v'+z'\big\rangle\\
&&\le (\|e\|_{L^4(D)}+\|\delta\|_{L^4(D)})
(\|v\|_{L^4(D)}+\|z\|_{L^4(D)}+\|v'\|_{L^4(D)}+\|z'\|_{L^4(D)})\|e\|_{\IH^1}.
\eeqas
We now use the following interpolation inequality 
\be
\|e\|_{L^4(D)}\le c_0\|e\|_{\IH^1}^{1/2}\|e\|_{\IH}^{1/2},
\label{eq:interpolationinequality}
\ee
which holds for all two dimensional domains $D$ with constant $c_0$ depending
only on $D$; see Flandoli \cite{Flandoli94}.
Using this we obtain
\begin{eqnarray*}
&&{1\over 2}{d\|e\|^2_{\IH}\over dt}+\nu\|e\|^2_{\IH^1}\le 
\\
&&\quad\quad c_1\Bigl(\|e\|_{\IH^1}^{3/2}\|e\|_{\IH}^{1/2}+\|\delta\|_{L^4(D)}\|e\|_{\IH^1}\Bigr) \cdot
\Bigl(\|v\|_{L^4(D)}+\|v'\|_{L^4(D)}+\|z\|_{L^4(D)}+\|z'\|_{L^4(D)}\Bigr)
\end{eqnarray*}
for a positive constant $c_1$. 
From the Young inequality, we have
\beqas
&&\|e\|^{3/2}_{\IH^1}\|e\|_{\IH}^{1/2}\Bigl(\|v\|_{L^4(D)}+\|v'\|_{L^4(D)}+\|z\|_{L^4(D)}+\|z'\|_{L^4(D)}\Bigr)\\
&&\qquad\le \frac34c_2^{4/3}\|e\|_{\IH^1}^2+\frac{1}{4c_2^4}\|e\|_{\IH}^2\Bigl(\|v\|_{L^4(D)}+\|v'\|_{L^4(D)}+\|z\|_{L^4(D)}+\|z'\|_{L^4(D)}\Bigr)^4
\eeqas
and
\beqas
&&\|\delta\|_{L^4(D)}\|e\|_{\IH^1}\Bigl(\|v\|_{L^4(D)}+\|v'\|_{L^4(D)}+\|z\|_{L^4(D)}+\|z'\|_{L^4(D)}\Bigr)\\
&&\qquad\le\frac{c_3^2}{2}\|e\|_{\IH^1}^2+\frac{1}{2c_3^2}\|\delta\|^2_{L^4(D)}\Bigl(\|v\|_{L^4(D)}+\|v'\|_{L^4(D)}+\|z\|_{L^4(D)}+\|z'\|_{L^4(D)}\Bigr)^2
\eeqas
for all positive constants $c_2$ and $c_3$. Choosing $c_2$ and $c_3$ so that
$c_1(3c_2^{4/3}/4+c_3^2/2)=\nu$, we deduce that there is a positive constant $c$ such that
\begin{equation}
{1\over 2}{d\|e\|_{\IH}^2\over dt}+\nu\|e\|_{\IH^1}^2\le
\nu\|e\|_{\IH^1}^2 +c\|e\|_{\IH}^2\cdot I_4
+c\|\delta\|_{L^4(D)}^2\cdot I_2.
\label{eq:eq}
\end{equation}
where we have defined
\begin{eqnarray*}
I_2&=\|v\|_{L^4(D)}^2+\|v'\|_{L^4(D)}^2+\|z\|_{L^4(D)}^2+\|z'\|_{L^4(D)}^2\\
I_4&=\|v\|_{L^4(D)}^4+\|v'\|_{L^4(D)}^4+\|z\|_{L^4(D)}^4+\|z'\|_{L^4(D)}^4.
\end{eqnarray*}
From Gronwall's inequality, we have
\begin{equation}
\|e(t)\|_{\IH}^2\le
c\int_0^t\bigg(e^{\int_s^t I_4(s')ds'} \bigg) 
\|\delta(s)\|_{L^4(D)}^2 I_2(s) ds.
\end{equation}
{Applying the interpolation inequality \eqref{eq:interpolationinequality} to $v(s';W)$,} 
we have that
\[
\int_0^T\|v(s';W)\|_{L^4(D)}^4 ds'\le c\sup_{s'\in(0,T)}\|v(s';W)\|_{\IH}^2\int_0^T\|v(s';W)\|_{\IH^1}^2ds',
\]
which is bounded uniformly when $W$ belongs to a bounded subset of $\cV$ due to \eqref{eq:bound}. {Using this estimate, and a similar estimate on $v'$,
together with  \eqref{eq:boundH12ofz} and Sobolev embedding of $\IH^{\frac12}$
into $L^4(D)$, we deduce that} 
\[
\|e(t)\|_{\IH}^2\le c\sup_{0\le s\le T}\|\delta(s)\|_{L^4(D)}^2.
\]

We now extend to include continuity with respect to the initial condition.
We show that $u(\cdot,t;u_0,W)$ is a continuous map from $\cH$ to $\IH$. For $W\in \cV$ and $u_0\in\IH$, we consider the following equation: 
\begin{equation}
\frac{dv}{dt}+A v+B(v+z,v+z)=0, \quad v(0)=u_0.
\label{eq:vu0}
\end{equation}
We denote the solution by $v(t)=v(t;u_0,W)$ to emphasize the dependence
on initial condition and forcing which is important here.
For $(u_0,W)\in\cH$ and $(u_0',W')\in\cH$, from \eqref{eq:eq} and Gronwall's inequality, we deduce that
\beqas
\|v(t;u_0,W)-v(t;u_0',W')\|^2_{\IH}\le
\|u_0-u_0'\|_{\IH}^2e^{\int_0^t I_4(s'))ds'}\\
\qquad\qquad+ c\int_0^t\bigg(e^{\int_s^t I_4(s')ds'}
\cdot \|z(s;W)-z(s;W')\|_{L^4(D)}^2. I_2(s)\bigg)ds.
\eeqas
We then deduce that
\beqas
\|v(t;u_0,W)-v(t;u_0',W')\|^2_{\IH}\le c\|u_0-u_0'\|^2_{\IH}+c\sup_{0\le s\le T}\|z(s;W)-z(s;W')\|_{L^4(D)}^2\\
\qquad\qquad\qquad\qquad\qquad\qquad\qquad\le c\|u_0-u_0'\|^2_{\IH}+c\sup_{t\in(0,T)}\|W(t)-W'(t)\|_{\IH^{1/2+\epsilon}}.
\eeqas
This gives the desired continuity of the forward map.
\eproof

\section{Bayesian Inverse Problems With Model Error}
\label{sec:B}

In this section we formulate the inverse
problem of determining the forcing to equation \eqref{eq:u}
from knowledge of the velocity field; more
specifically we formulate the Bayesian inverse problem of determining
the driving Brownian motion $W$ from noisy pointwise observations of the velocity field.
Here we consider the initial condition to be fixed and hence
denote the solution of \eqref{eq:u} by $u(t;W)$; extension to
the inverse problem for the pair $(u_0,W)$ is given in the
following section.

We set-up the likelihood in subsection \ref{ssec:like}.
Then, in subsection \ref{ssec:prior}, we describe
the prior on the forcing which is a Gaussian white-in-time
process with spatial correlations, and hence a spatially correlated
Brownian motion prior on $W$. This leads, in
subsection \ref{ssec:post}, to a well-defined
posterior distribution, absolutely continuous with
respect to the prior, and Lipschitz in the Hellinger
metric with respect to the data. {To prove these results
we employ the framework for Bayesian inverse problems developed in Cotter et al. \cite{CDRS} and Stuart \cite{stuart2010inverse}. In particular, Corollary 2.1 of \cite{CDRS} and Theorem 6.31 of \cite{stuart2010inverse} show that, in order to demonstrate the absolute continuity 
of the posterior measure with respect to the prior, it suffices to show that the mapping $\cG$ in \eqref{eq:delta} is continuous with respect to the topology of $\IX$ and to choose a prior with full mass on $\IX$.
Furthermore we then employ the proofs of Theorem 2.5 of \cite{CDRS} and Theorem 4.2 of \cite{stuart2010inverse} to show the well-posedness of the posterior measure;
indeed we show that the posterior is Lipschitz with respect to data, in the
Hellinger metric.}

\subsection{Likelihood}
\label{ssec:like}

Fix a set of times $t_j\in (0,T)$,
$j=1,\cdots, J.$ 
Let $\ell$ be a collection of $K$ bounded linear functionals
on $\IH$. We assume that we observe, for each $j$,
$\ell\bigl(u(\cdot,t_j;W)\bigr)$ 
plus a  draw from a centered $K$ dimensional Gaussian noise $\vartheta_j$ 
so that
\be
\delta_j=\ell\bigl(u(\cdot,t_j;W)\bigr)+\vartheta_j,
\label{eq:detlax}
\ee
is known to us. Concatenating the data we obtain
\be
\delta=\cG(W)+\vartheta
\label{eq:delta}
\ee
where $\delta,\vartheta \in \bbR^{JK}$ and $\cG: \cV \to \bbR^{JK}.$
The observational noise $\vartheta$ 
is a draw from the $JK$ dimensional Gaussian random variable with the covariance matrix $\Sigma$.

{
In the following we will define a prior measure $\rho$ on
$W$ and then determine the conditional probability measure 
$\rho^\delta=\IP(W|\delta)$. We will then show that 
$\rho^{\delta}$ is absolutely continuous with respect to $\rho$ 
and that the Radon-Nikodym derivative between the measures
is given by}
\be
{d\rho^\delta\over d\rho}\propto\exp\Bigl(-\Phi(W;\delta)\Bigr),
\label{eq:RN}
\ee
where 
\be
\Phi(W;\delta)={1\over 2}\Big|\Sigma^{-\frac12}\bigl(\delta-\cG(W)\big)\Big|^2.
\label{eq:Phi}
\ee
The right hand side of \eqref{eq:RN} is the likelihood of the data $\delta$.

\subsection{Prior}
\label{ssec:prior}

We construct our prior on the time-integral of the forcing, namely
$W$. 
Let $Q$ be a linear operator from the Hilbert space $\IH^{\frac12+\epsilon}$ into itself with eigenvectors $e_k$ and eigenvalues $\sigma_k^2$ for $k=1,2,\ldots$. We make the following assumption
\begin{assumption}\label{ass:noise}
There is an $\epsilon>0$ such that
the coefficients $\{\sigma_k\}$ satisfy
\[
\sum_{k=1}^\infty \sigma_k^2\lambda_k^{1/2+\epsilon}<\infty
\]
where $\lambda_k$ are the eigenvalues of the operator $A$ defined in Section \ref{ssec:over}.
\end{assumption}
As 
\[
\sum_{k=1}^\infty\langle Qe_k,e_k\rangle_{\IH^{\frac12+\epsilon}}=\sum_{k=1}^\infty\sigma_k^2\lambda_k^{\frac12+\epsilon}<\infty,
\]
$Q$ is a trace class operator in $\IH^{\frac12+\epsilon}$.  

We assume that our prior is the $Q$-Wiener process $W$ with values in $\IH^{\frac12+\epsilon}$ where $W(s_1)-W(s_2)$ is Gaussian in $\IH^{\frac12+\epsilon}$ with covariance $(s_1-s_2)Q$ and mean 0. This process can be written as 
\be
W(t)=\sum_{k=1}^\infty \sigma_k e_k\omega_k(t),
\label{eq:w}
\ee
where $\omega_k(t)$ are pair-wise independent Brownian motions (see Da Prato and Zabczyk \cite{DPZ}, Proposition 4.1) and where the convergence of the infinite series is in the mean square norm with respect to the probability measure of the probability space that generates the randomness of $W$.  
We define by $\rho$ the measure generated by this $Q$-Wiener process on $\cV$.

\begin{remark}
{We have constructed the solution to \eqref{eq:u} for each deterministic continuous function $W\in \IX$. As we equip $\IX$ with the prior probability measure $\rho$, we wish to employ the results from \cite{Flandoli94} concerning
the solution of \eqref{eq:u} when $W$ is considered as a Brownian motion obtaining values in $\IX$. However,  the solution of \eqref{eq:u} is constructed in a 
slightly different way in \cite{Flandoli94} from that used in the preceding
developments.  
We therefore show that} under Assumption \ref{ass:noise}, $\rho$ almost surely, solution $u$ of \eqref{eq:weaksln} defined in \eqref{eq:defu} for each individual function $W$ equals the unique progressively measurable solution in $C([0,T];\IH)\bigcap L^2([0,T];\IH^1)$ constructed in Flandoli \cite{Flandoli94} when the noise $W$ is sufficiently spatially regular. This allows us to employ the existence of the second moment of $\|u(\cdot,t;W)\|_{\IH}^2$, i.e.the the finiteness of the energy $\IE^{\rho}[\|u(\cdot,t;W)\|_{\IH}^2]$, established in Mattingly \cite{Mattingly99}, which we need later.

For the infinite dimensional Brownian motion $W$ defined in \eqref{eq:w} where 
\[
\sum_{k=1}^\infty \lambda_k^{2\beta_0-1/2}\sigma_k^2<\infty,
\]
for some $\beta_0>0$, {where we employ the same notation 
as in \cite{Flandoli94} for ease of exposition.} 
Flandoli \cite{Flandoli94} employs the Ornstein-Uhlenbeck process
\be
z_\alpha(t)=\int_{-\infty}^te^{-(\nu A+\alpha)(t-s)}dW(s)
\label{eq:OU}
\ee
which, considered as the stochastic process, is a solution of the Ornstein-Uhlenbeck equation
\[
dz_\alpha(t)+Az_\alpha(t)dt+\alpha z_\alpha(t)dt=dW(t),
\]
where $\alpha$ is a constant, {in order to define
a solution of \eqref{eq:weaksln}. Note that if $\beta_0>\frac12$
then Assumption \ref{ass:noise} is satisfied.}
With respect to the probability space $(\Omega,\cF_t,\IP)$, the expectation $\IE\|z_\alpha(t)\|^2_{\IH^{1/2+2\beta}}$  is finite for $\beta<\beta_0$. Thus almost surely with respect to $(\Omega,\cF_t,\IP)$, $z_\alpha(t)$ is sufficiently regular so that problem \eqref{eq:v} with the initial condition $v(0;W)=u_0-z_\alpha(0)$ is well posed. The stochastic solution to the problem \eqref{eq:u} is defined as
\be
u(\cdot,t;W)=z_\alpha(t;W)+v(t;W)
\label{eq:ualpha}
\ee
which is shown to be independent of $\alpha$ in \cite{Flandoli94}.
When $\beta_0>\frac12$, $\IE\|z_\alpha(t)\|^2_{\IH^1}$ is finite so $u(\cdot,t;W)\in C([0,T];\IH)\bigcap L^2([0,T];\IH^1)$. Flandoli \cite{Flandoli94} leaves open the question of the uniqueness of a generalized solution to \eqref{eq:weaksln} in $C([0,T];\IH)\bigcap L^2([0,T];\IH^{1/2})$. However, there is a unique solution in $C([0,T];\IH)\bigcap L^2([0,T];\IH^1)$.

Almost surely with respect to the probability measure $\rho$, solution $u$ of \eqref{eq:weaksln} constructed in \eqref{eq:defu} equals the solution constructed by Flandoli \cite{Flandoli94} in \eqref{eq:ualpha}.
To see this, note that the stochastic integral
\be
\int_0^t e^{-\nu A(t-s)}dW(s)
\label{eq:0tOU}
\ee
can be written in the integration by parts form \eqref{eq:zz}. Therefore, with respect to $\rho$, 
\beqas
\IE^{\rho}[\|z(t)\|_{L^2(0,T;\IH^1)}^2]&=&\int_0^T\IE^\rho[\sum_{k=1}^\infty\lambda_k\sigma_k^2\Bigg(\int_0^te^{-\nu\lambda_k(t-s)}d\omega_k(s)\Bigg)^2]dt\\
&=&\int_0^T\Bigg(\sum_{k=1}^\infty\lambda_k\sigma_k^2\int_0^te^{-2\nu\lambda_k(t-s)}ds\Bigg)dt={1\over2\nu}\int_0^T\sum_{k=1}^\infty \sigma_k^2(1-e^{-2\nu\lambda_kt})dt
\eeqas
which is finite. Therefore $\rho$ almost surely, $z(t)\in L^2(0,T;\IH^1)$. Thus $u(t;W)\in C(0,T;\IH)\bigcap L^2(0,T;\IH^1)$. We can then argue that $\rho$ almost surely, the solution $u$ constructed in \eqref{eq:defu} equals Flandoli's solution in \eqref{eq:ualpha} which we denote by $u_\alpha$ (even though it does not depend on $\alpha$) as follows. 
As $u_\alpha\in C([0,T];H)\bigcap L^2([0,T];\IH^1)$, $v_\alpha(t;W)=u_\alpha(t;W)-z(t;W)\in 
C([0,T];H)\bigcap L^2([0,T];\IH^1)$ and satisfies \eqref{eq:v}. As for each $W$, \eqref{eq:v} has a unique solution in $C([0,T];H)\bigcap L^2([0,T];\IH^1)$, so $v_\alpha(t;W)=v(t;W)$. 
Thus almost surely, the Flandoli \cite{Flandoli94} solution equals the solution $u$ in \eqref{eq:defu}.   
This is also the argument to show that \eqref{eq:u} has a unique solution in $C([0,T];\IH)\bigcap L^2([0,T];\IH^1)$.


\end{remark}

\subsection{Posterior}
\label{ssec:post}

{
\begin{theorem}\label{thm:posterior}
The conditional measure $\IP(W|\delta)=\rho^\delta$ is absolutely continuous with respect to the prior measure $\rho$ with the Radon-Nikodym derivative being given by \eqref{eq:RN}.
Furthermore, there is a constant $c$ so that
\[
d_{\rm Hell}(\rho^\delta,\rho^{\delta'})\le c|\delta-\delta'|.
\]
\end{theorem}
}

\bproof
{Note that $\rho(\IX)=1.$ It follows
from Corollary 2.1 of Cotter et al. \cite{CDRS} and Theorem 6.31 of Stuart \cite{stuart2010inverse} that, in order to demonstrate that $\rho^{\delta} \ll \rho$,
it suffices to show that the mapping $\cG: \cV \to \bbR^{JK}$ is continuous;
then the Randon-Nikodym derivative \eqref{eq:RN} defines the density of
$\rho^{\delta}$ with respect to $\rho.$
As $\ell$ is a collection of bounded continuous linear functionals on $\IH$, the continuity of $\cG$ with respect to the topology of $\IX$ follows from Theorem \ref{t:first}.

We now turn to the Lipschitz continuity of the posterior in the Hellinger metric.
The method of proof is very similar to that developed in the 
proofs of Theorem 2.5 in \cite{CDRS} and Theorem 4.2 in \cite{stuart2010inverse}. 
We define
\[
Z(\delta):=\int_{\cV}\exp(-\Phi(W;\delta))d\rho(W).
\]
Mattingly \cite{Mattingly99} shows that for each $t$, the second moment $\IE^\rho(\|u(\cdot,t;W)\|_{\IH}^2)$ is finite. Fixing a large constant 
$M$, the $\rho$ measure of the set of paths $W$ such  that $\max_{j=1,\ldots,J}\|u(\cdot,t_j;W)\|_{\IH}\le M$ is 
larger than $1-cJ/M^2>1/2$. For those paths $W$ in this set we have,
\[
\Phi(W;\delta)\le c(|\delta|+M).
\]
From this, we deduce that $Z(\delta)>0$. 
Next, we have that
\beqas
|Z(\delta)-Z(\delta')|&\le& \int_{\cV}|\Phi(W;\delta)-\Phi(W;\delta')|d\rho(W)\\
&\le& c\int_{\cV}(|\delta|+|\delta'|+2\sum_{j=1}^J|\ell(u(t_j;W))|_{\IR^K})|\delta-\delta'|d\rho(W)\\
&\le& c|\delta-\delta'|.
\eeqas
We then have
\beqas
2d_{\rm Hell}(\rho^\delta,\rho^{\delta'})^2&\le& \int_{\cV}\left(Z(\delta)^{-1/2}\exp(-{1\over 2}\Phi(W;\delta))-Z(\delta')^{-1/2}\exp(-{1\over 2}\Phi(W';\delta'))\right)^2d\rho(W)\\
&\le& I_1+I_2,
\eeqas
where
\[
I_1={2\over Z(\delta)}\int_{\cV}\left(\exp(-{1\over 2}\Phi(W;\delta))-\exp(-{1\over 2}\Phi(W;\delta'))\right)^2d\rho(W),
\]
and 
\[
I_2=2|Z(\delta)^{-1/2}-Z(\delta')^{-1/2}|^2\int_{\cV}\exp(-\Phi(W;\delta'))d\rho(W).
\]
Using the facts that
\beqas
\Big|\exp(-{1\over 2}\Phi(W;\delta))-\exp(-{1\over 2}\Phi(W;\delta'))\Big|\le
\frac12|\Phi(W;\delta)-\Phi(W;\delta')|
\eeqas
and that $Z(\delta)>0$, we deduce that
\beqas
I_1&\le& c\int_\IX|\Phi(W;\delta)-\Phi(W;\delta')|^2d\rho(W)\\
&\le& c\int_\cV(|\delta|+|\delta'|+2\sum_{j=1}^J|\ell(u(t_j;W))|_{\IR^K})^2|\delta-\delta'|^2d\rho(W)\\
&\le& c|\delta-\delta'|^2.
\eeqas
Furthermore,
\[
|Z(\delta)^{-1/2}-Z(\delta')^{-1/2}|^2\le c\max(Z(\delta)^{-3},Z(\delta')^{-3})|Z(\delta)-Z(\delta')|^2.
\]
From these inequalities it follows
that $d_{\rm Hell}(\rho^\delta,\rho^{\delta'})\le c|\delta-\delta'|.$
\eproof

\section{Inferring The Initial Condition}
\label{sec:C}

In the previous section we discussed the problem of inferring
the forcing from the velocity field. In practical applications it
is also of interest to infer the initial condition, which corresponds
to a Bayesian interpretation of 4DVAR, or the initial
condition and the forcing, which corresponds
to a Bayesian interpretation of weak constraint 4DVAR.
Thus we consider {the} Bayesian inverse problem for inferring the initial condition $u_0$ and the 
white noise forcing determined by the Brownian driver $W$. 
{Including the initial condition does not add any further
technical difficulties as the dependence on the pathspace valued
forcing is more subtle than the dependence on initial condition,
and this dependence on the forcing is dealt with in the previous
section. As a consequence we do not provide full details.}

Let $\varrho$ be a Gaussian measure on the space $\IH$ and let $\mu=\varrho\otimes\rho$ be the prior probability measure on the space $\cH=\IH\times\cV$.  We denote the solution $u$ of \eqref{eq:u} by $u(x,t;u_0,W)$. 

We outline what is required to extend the analysis of the previous two sections
to the case of inferring both initial condition and driving Brownian motion.
We simplify the presentation by assuming observation at only one time $t_0>0$ 
although this is easily relaxed. 
Given that at $t_0\in (0,T)$, the noisy observation $\delta$ of $\ell(u(\cdot,t_0;u_0,W)$ is given by 
\be
\delta=\ell(u(\cdot,t_0;u_0,W))+\vartheta
\label{eq:deltau0}
\ee
where $\vartheta\sim N(0,\Sigma)$. Letting 
\be
\Phi(u_0,W;\delta)={1\over2}|\delta-\ell(u(\cdot,t_0;u_0,W))|_\Sigma^2,
\ee
we aim to show that the conditional probability $\mu^\delta=\IP(u_0,W|\delta)$ is given by
\be
{d\mu^\delta\over d\mu}\propto\exp(-\Phi(u_0,W;\delta)).
\label{eq:RNu0}
\ee
We have the following result.

\begin{theorem}
The conditional probability measure $\mu^\delta=\IP(u_0,W|\delta)$ is absolutely continuous with respect to the prior probability measure $\mu$ with the Radon-Nikodym derivative give by \eqref{eq:RNu0}. Further, there is a constant $c$ such that 
\[
d_{Hell}(\mu^\delta,\mu^{\delta'})\le c|\delta-\delta'|.
\] 
\end{theorem}

\bproof
{To establish the absolute continuity of posterior with respect to
prior, together with the formula for the Radon-Nikodym derivative,
the key issue is establishing continuity of the forward
map with respect to initial condition and driving Brownian motion
This is established in Theorem \ref{t:first} and since $\mu(\cH)=1$
the first part of the theorem follows.}

For the Lipschitz dependency of the Hellinger distance of $\mu^\delta$ on $\delta$, we use the result of 
 Mattingly \cite{Mattingly99} which shows that, for each initial condition $u_0$, 
\[
\IE^\rho(\|u(t;u_0,W)\|_{\IH}^2)\le {{\mathcal E}_0\over2\nu\lambda_1}+e^{-2\nu\lambda_1 t}(\|u_0\|_{\IH}^2-{{\mathcal E}_0\over 2\nu\lambda_1}),
\]
where ${\mathcal E}_0=\sum_{k=1}^\infty\sigma_k^2$. Therefore $\IE^\mu(\|u(t;u_0,W)\|_{\IH}^2)$ is bounded. {This enables us to establish positivity of the
normalization constants and the remainder of the proof follows that
given in Theorem \ref{thm:posterior}.}
\eproof

\section{Numerical Results}
\label{sec:D}

{The purpose of this section is twofold:
firstly to demonstrate that the
Bayesian formulation of the inverse problem described in this
paper forms the basis for practical numerical inversion;
and secondly to study some properties of the posterior
distribution on the white noise forcing, given observations of
linear functionals of the velocity field.

The numerical results move outside the strict remit
of our theory in two directions.
Firstly we work with periodic boundary conditions; this
makes the computations fast, but simultaneously demonstrates
the fact that the theory is readily extended from Dirichlet
to other boundary conditions. Secondly we consider both 
(i) pointwise observations of the entire velocity field and
(ii) observations found from the projection onto the lowest eigenfunctions 
of $A$ noting that the second form of observations are
bounded linear functionals on  $\IH$, as required by our theory,
whilst the first form of observations are not.

To extend our theory to periodic boundary
conditions requires generalization of the Flandoli \cite{Flandoli94}
theory from the Dirichlet to the periodic setting, which is not a 
technically challenging generalization. However consideration
of pointwise observation functionals requires the proof
of continuity of $u(t;\cdot, \cdot)$ as a mapping from $\cH$
into $\IH^s$ spaces for $s$ sufficiently large. Extension of
the theory to include pointwise
observation functionals would thus involve significant technical
challenges, in particular to derive smoothing estimates for the
semigroup underlying the Flandoli solution concept. 
Our numerical results will show that the posterior distribution for (ii)
differs very little from that for (i), which is an interesting fact in its
own right. 
}

In subsection \ref{num:fwd} we describe the numerical method used
for the forward problem. In subsection \ref{num:ips} we describe
the inverse problem and the Metropolis-Hastings MCMC method used to
probe the posterior. 
Subsection \ref{num:resdis} describes the  numerical results.

\subsection{Forward Problem: Numerical Discretization}
\label{num:fwd}

All our numerical results are computed using a viscosity
of $\nu=0.1$ and on the periodic domain.
We work on the time interval $t \in [0,0.1].$
We use $M=32^2$ divergence free Fourier basis functions for
a spectral Galerkin 
spatial approximation, 
and employ a time-step $\delta t = 0.01$ in 
a Taylor time-approximation \cite{jentzen2010taylor}.
The number of basis functions and time-step lead to
a fully-resolved numerical simulation at this value of $\nu.$

\subsection{Inverse Problem: Metropolis Hastings MCMC}
\label{num:ips}

Recall the Stokes' operator $A$.
We consider the inverse problem of finding the driving Brownian
motion. {As a prior we take a centered Brownian motion in time
with spatial covariance $\pi^4 A^{-2}$;} thus the space-time
covariance of the process is
$C_0 := \pi^4 A^{-2} \otimes (-\triangle_t)^{-1}$, where $\triangle_t$ is
the Laplacian in time with fixed homogeneous Dirichlet condition at
$t=0$ and homogeneous Neumann condition at $t=T$.
It is straightforward to draw samples from this Gaussian measure,
using the fact that $A$ is diagonalized in the spectral basis.
Note that if $W \sim \rho$, then $W \in C(0,T;\bbH^s)$ almost surely
for all $s<1$; in particular $W \in \bbX$.  Thus $\rho(\bbX)=1$ as
required. The likelihood is defined (i) by making observations
of the velocity field at every point on the $32^2$ grid implied
by the spectral method, at every time $t=n\delta t$, $n=1,\cdots 10$,
or (ii) by making observations of the projection onto eigenfunctions 
$\{\phi_k\}_{|k|<4}$ of $A$.
The observational noise standard deviation is taken to be
$\gamma = 1.6$ and all observational noises are uncorrelated.}

To sample from the posterior distribution we employ
a Metropolis-Hastings MCMC method. Furthermore, to ensure
mesh-independent convergence properties, we use a method
which is well-defined in function space \cite{CRSW12}. 
Metropolis-Hastings methods proceed by constructing
a Markov kernel $\cP$ which satisfies {\it detailed balance} with
respect to the measure $\rho^{\delta}$ which we wish to sample: 
\begin{equation}
\rho^{\delta}(du) \cP(u,dv) = \rho^{\delta}(dv) \cP(v,du), \quad \forall ~ u,v \in \bbX.
\label{db}
\end{equation}
Integrating with respect to $u$, one can see that 
detailed balance implies $\rho^{\delta}\cP = \rho^{\delta}$.
Metropolis-Hastings methods \cite{Has70,Tie} 
prescribe an accept-reject move based
on proposals from another Markov kernel $\cQ,$ in order 
to define a kernel $\cP$ which satisfies detailed balance.
If we define the measures 
\begin{eqnarray}
\begin{array}{lll}
\nu(du,dv) &=& \cQ(u,dv) \rho^\delta(du) \propto \cQ(u,dv) \exp\Bigl(-\Phi(u;\delta)\Bigr)\rho(du) \\
\nu^\perp(du,dv) &=& \cQ(v,du) \rho^\delta(dv) \propto \cQ(v,du) \exp\Bigl(-\Phi(v;\delta)\Bigr)\rho(dv).
\end{array}
\end{eqnarray}
then, provided $\nu^\perp \ll \nu$, the Metropolis-Hastings method
is defined as follows.  Given current state $u_n$,
a proposal is drawn $u^* \sim \cQ(u_n,\cdot)$, 
and then accepted with probability
\begin{equation} 
\alpha(u_n,u^*) = {\rm min} 
\left \{1, 
\frac{d\nu^\perp}{d\nu}(u_n,u^*)
 \right \}.
\label{acceptance1} 
\end{equation}
The resulting chain is denoted by $\cP.$ 
If the proposal $\cQ$ preserves the prior, so that $\rho \cQ=\rho$,
then a short calculation reveals that
\begin{equation} 
\alpha(u_n,u^*) = {\rm min} 
\left \{1, 
\exp\Bigl(\Phi(u_n;\delta)-\Phi(u^*;\delta)\Bigr)
 \right \};
\label{acceptance} 
\end{equation}
thus the acceptance probability is determined by the
change in the likelihood in moving from current to
proposed state.
We use the following pCN proposal \cite{CRSW12}
which is reversible with respect to the Gaussian
prior $N(0,C_0)$:
\begin{equation}
\cQ(u_n,\cdot) = N\bigl(\sqrt{1-\beta^2} u_n, \beta^2 C_0\bigr).
\label{pcnprop}
\end{equation}
This hence results in the acceptance probability \eqref{acceptance}.
Variants on this algorithm, which propose differently in
different Fourier components, are described in 
\cite{law2012proposals}, and can make substantial speedups
in the Markov chain convergence. However for the examples
considered here the basic form of the method suffices.

\subsection{Results and Discussion}
\label{num:resdis}

The true driving Brownian motion $W^\dagger$, 
underlying the data in the likelihood,
is constructed as a draw from the prior
$\rho$. We then compute the corresponding 
true trajectory $u^\dagger(t)=u(t;W^\dagger)$.  
We use the pCN scheme
\eqref{acceptance},\eqref{pcnprop} to sample $W$ from 
the posterior distribution $\rho^\delta$.
It is important to appreciate that the object of
interest here is the posterior distribution on $W$ itself
which provides estimates of the forcing, given the noisy
observations of the velocity field. This posterior distribution
is not necessarily close to a Dirac measure on the truth;
in fact we will show that some parameters required to
define $W$ are recovered accurately whilst others are not.

\begin{figure}
 \includegraphics[width=0.5\textwidth]{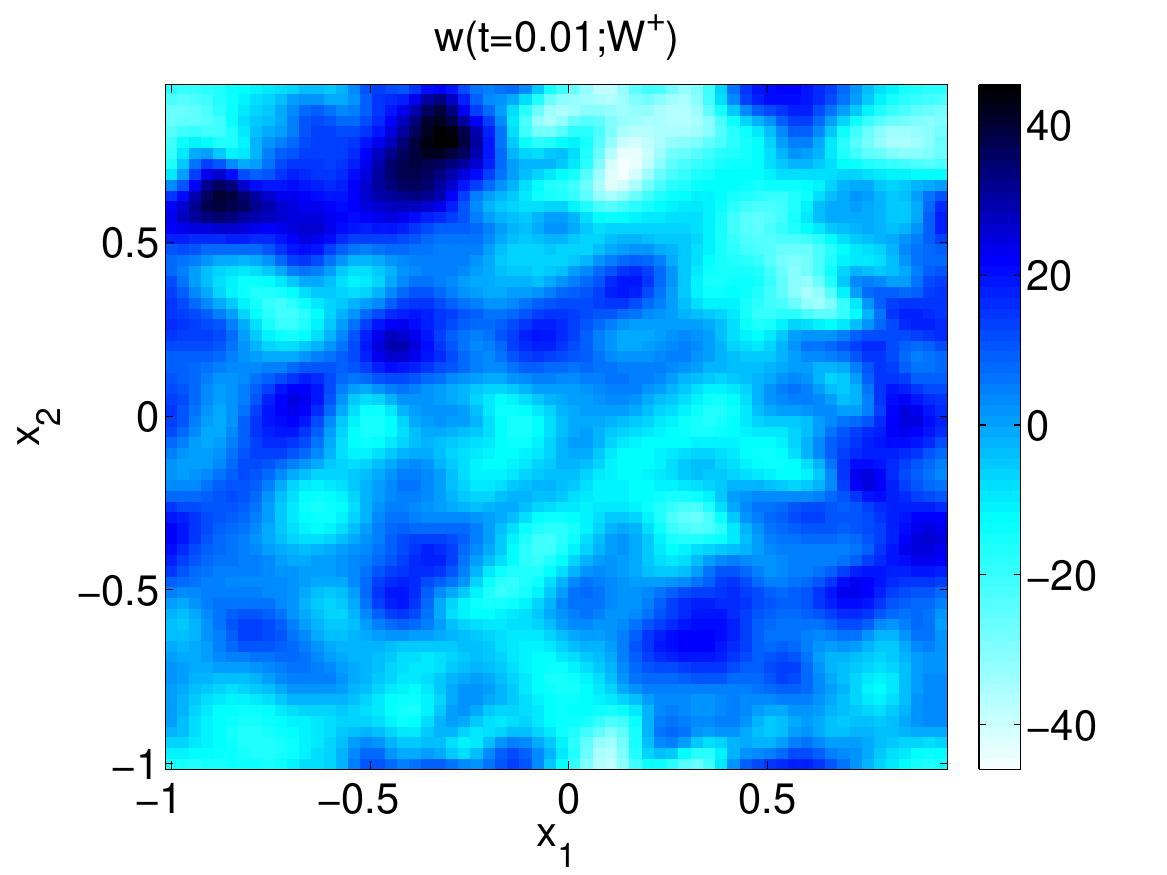}
 \includegraphics[width=0.5\textwidth]{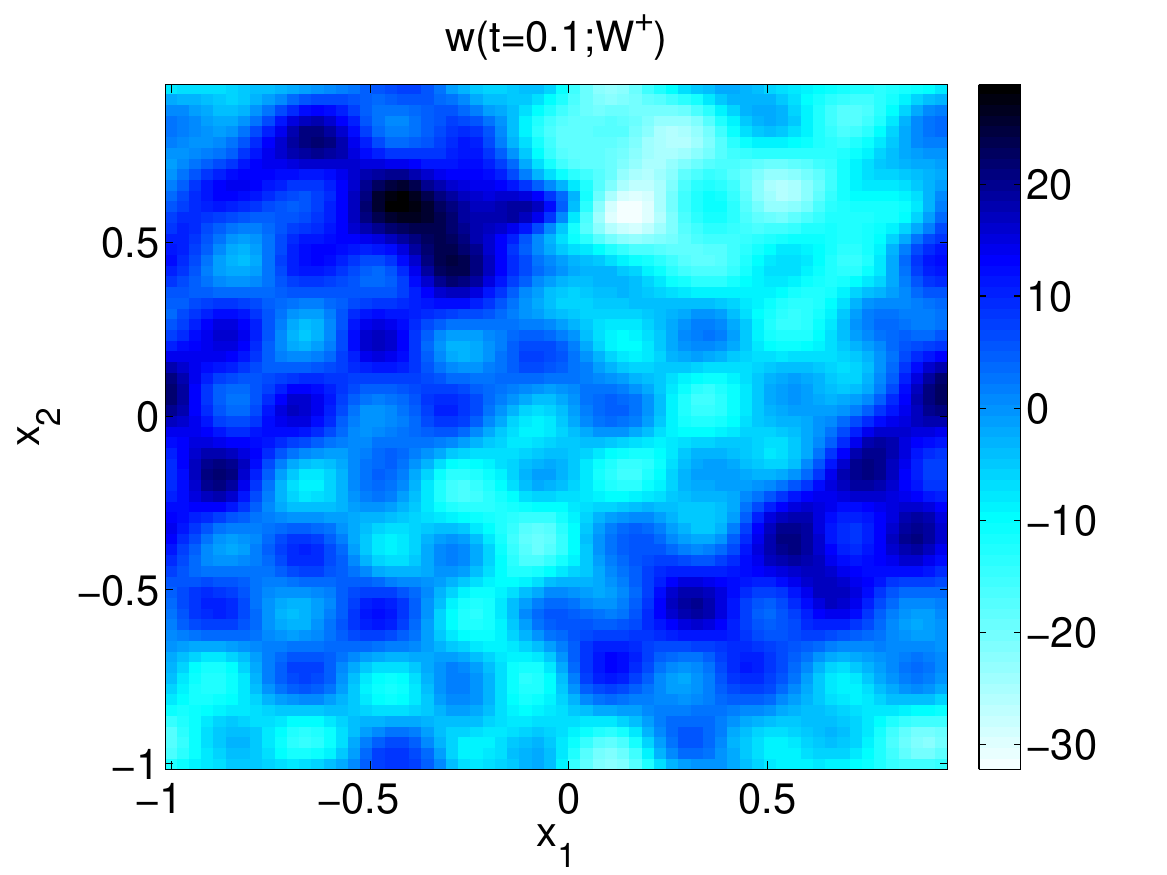}
 \includegraphics[width=0.5\textwidth]{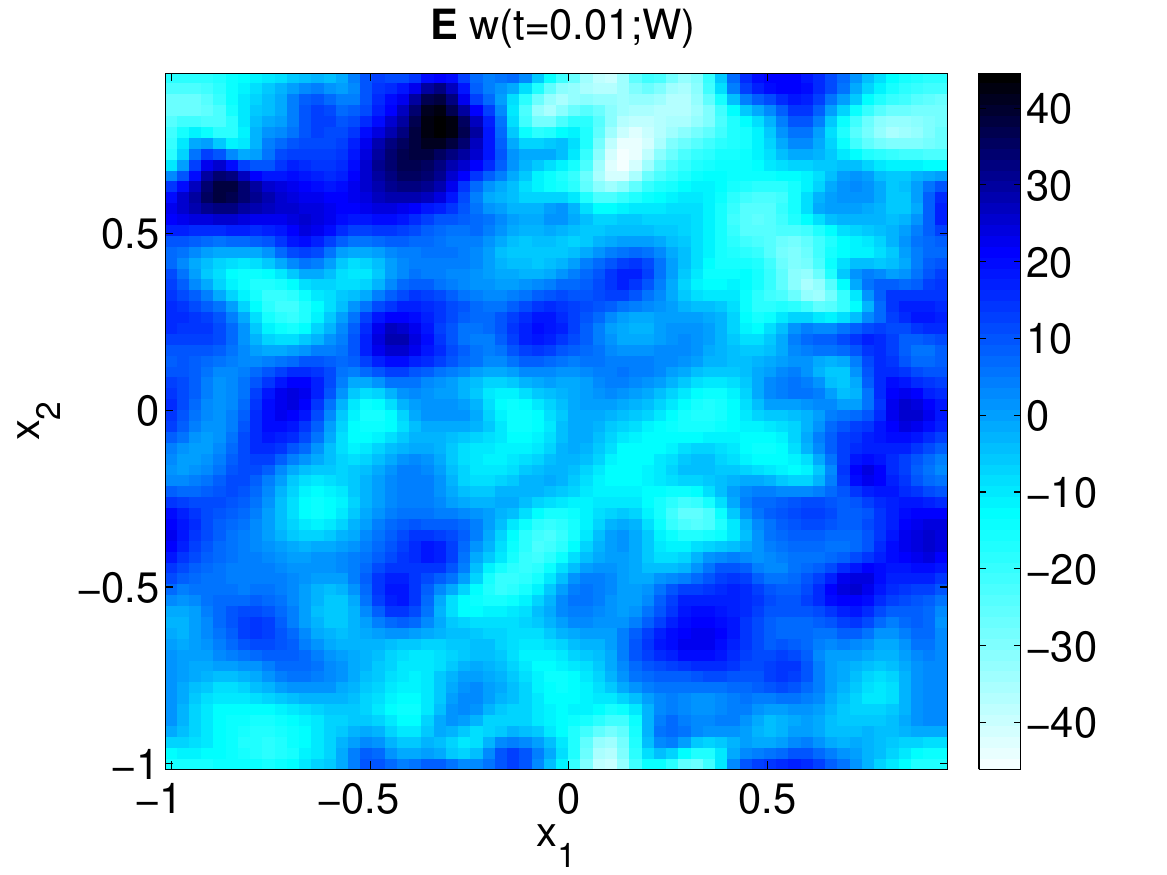}
 \includegraphics[width=0.5\textwidth]{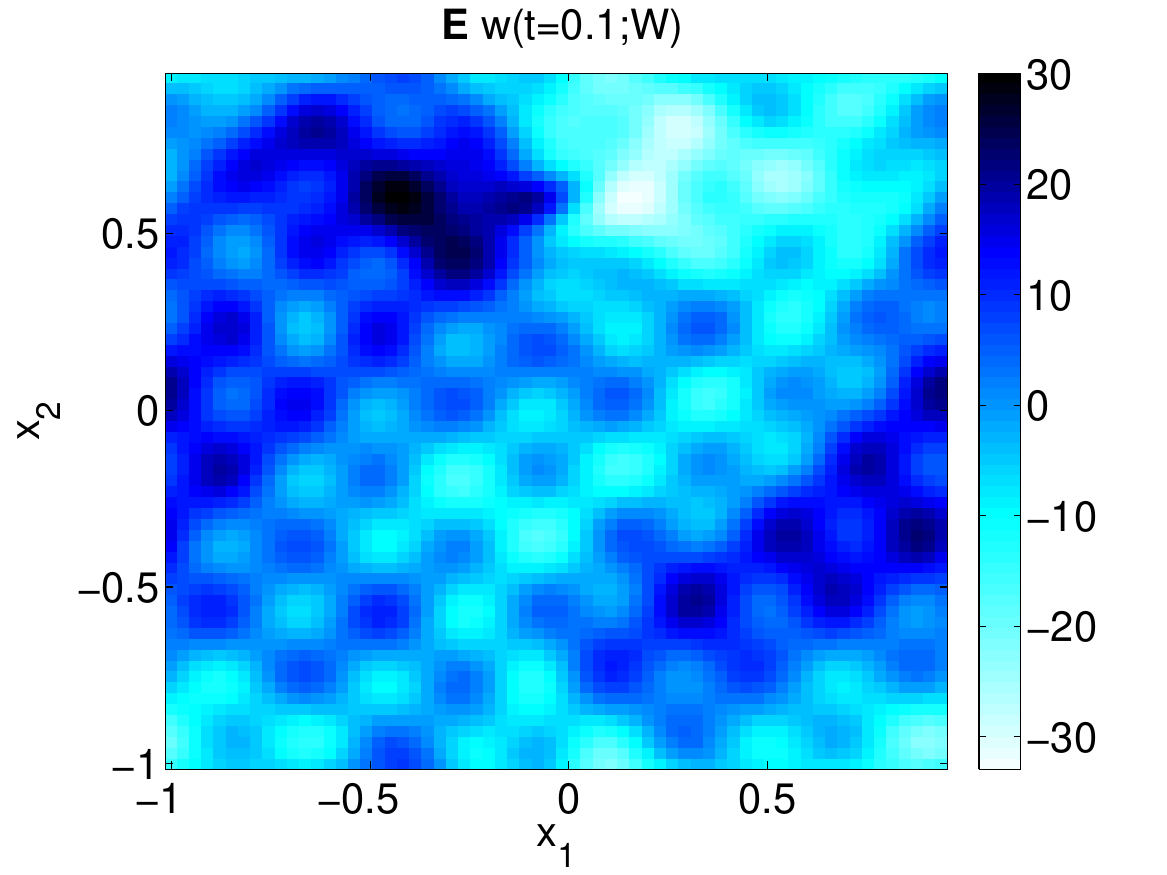}
 \includegraphics[width=0.5\textwidth]{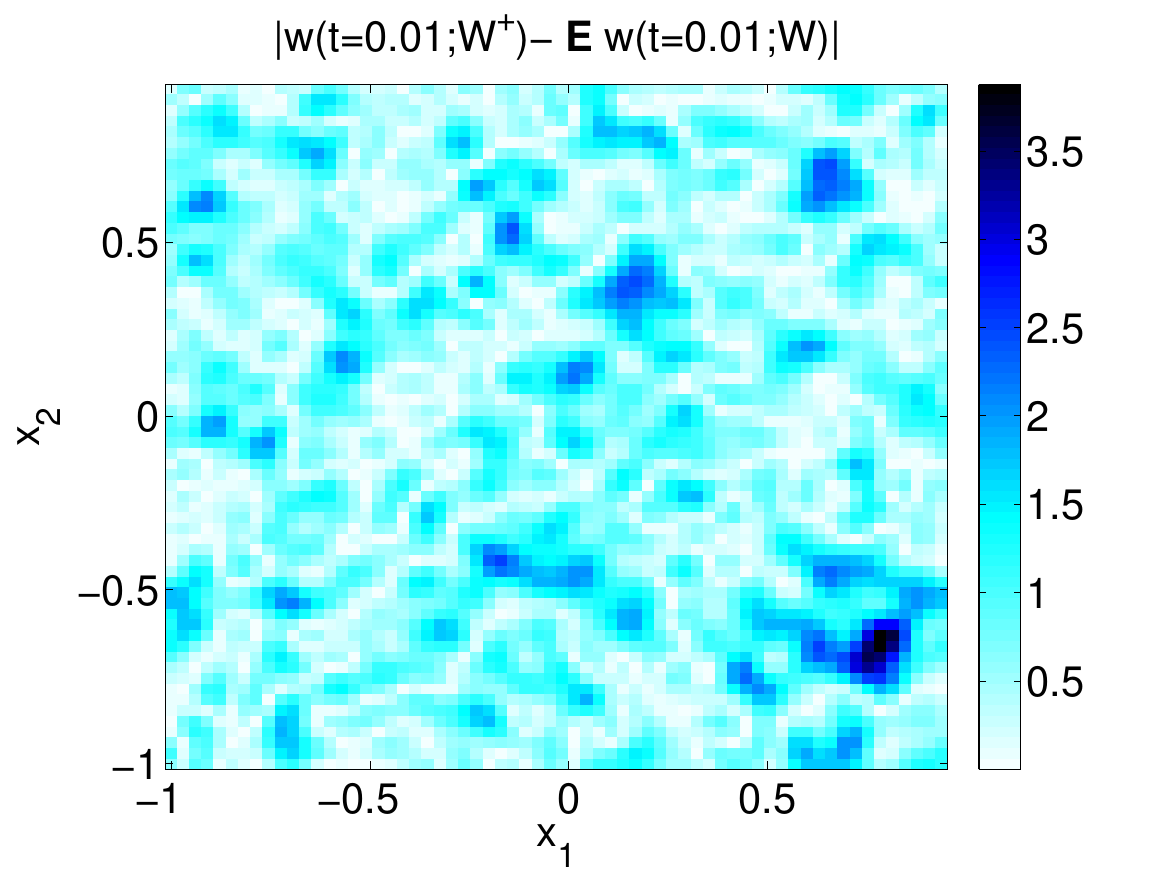}
 \includegraphics[width=0.5\textwidth]{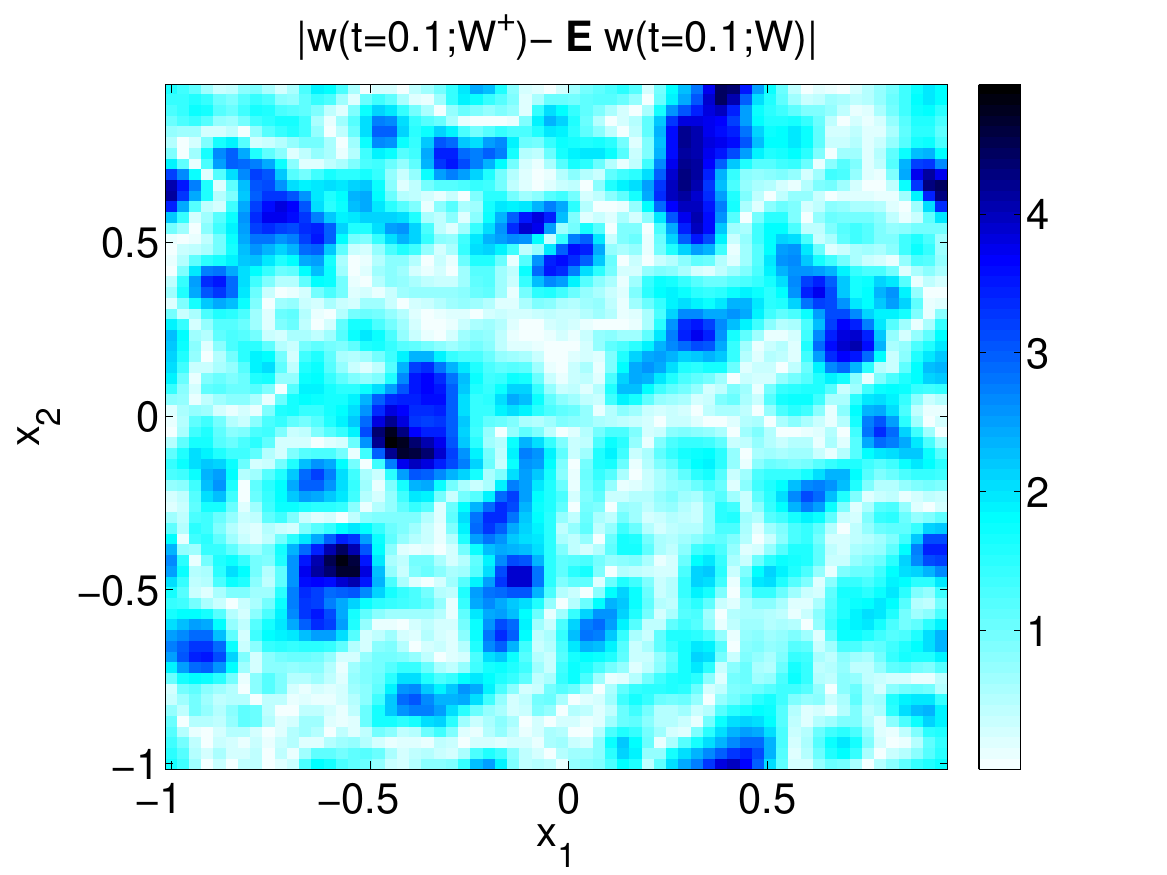} 
 \caption{Full-field point wise observations. 
 The truth (top), expected value (middle), and
 absolute distance between them (bottom) of the vorticity
 $w(t;W)$, for $t=0.01$ 
 (left, relative $L^2$ error $e=0.0044$) 
 and $t=0.1$ (right, $e=0.0244$).}
\label{prof}
\end{figure}

\begin{figure}
 \includegraphics[width=0.5\textwidth]{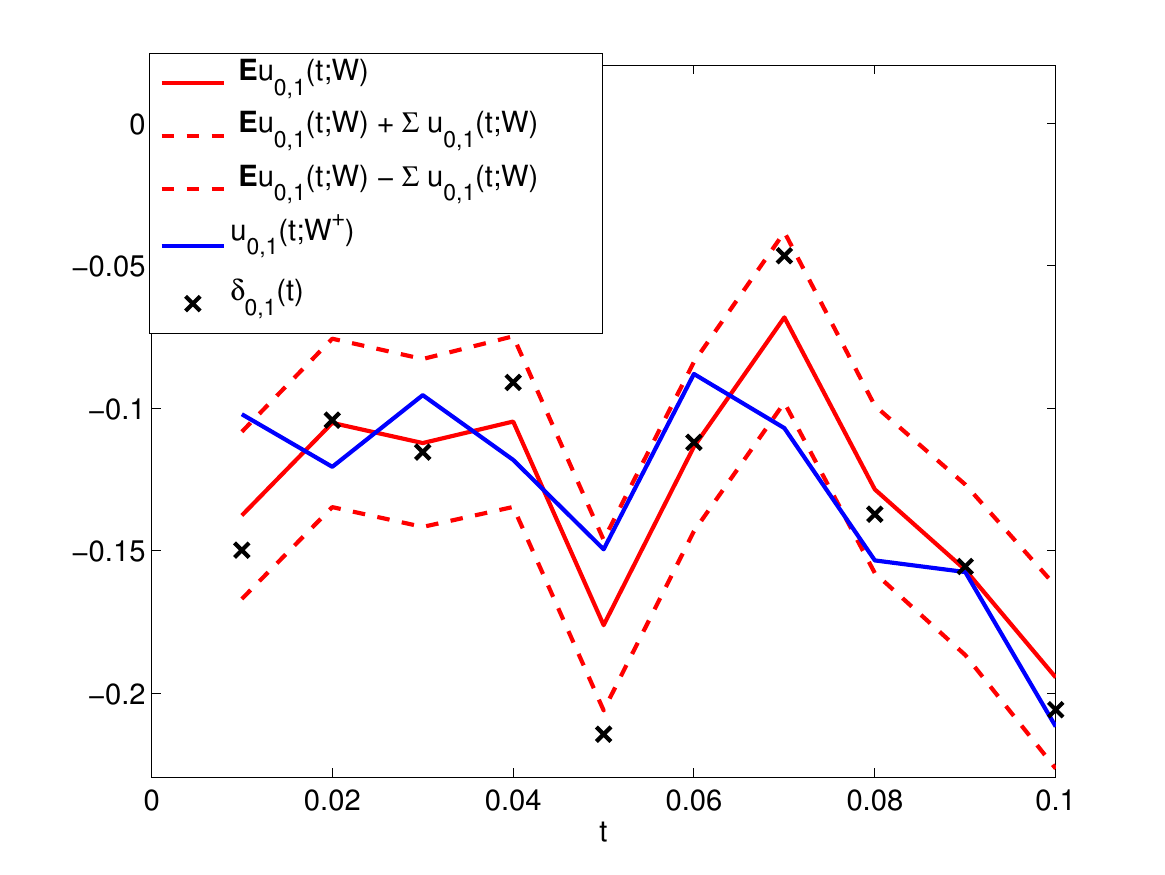}
 \includegraphics[width=0.5\textwidth]{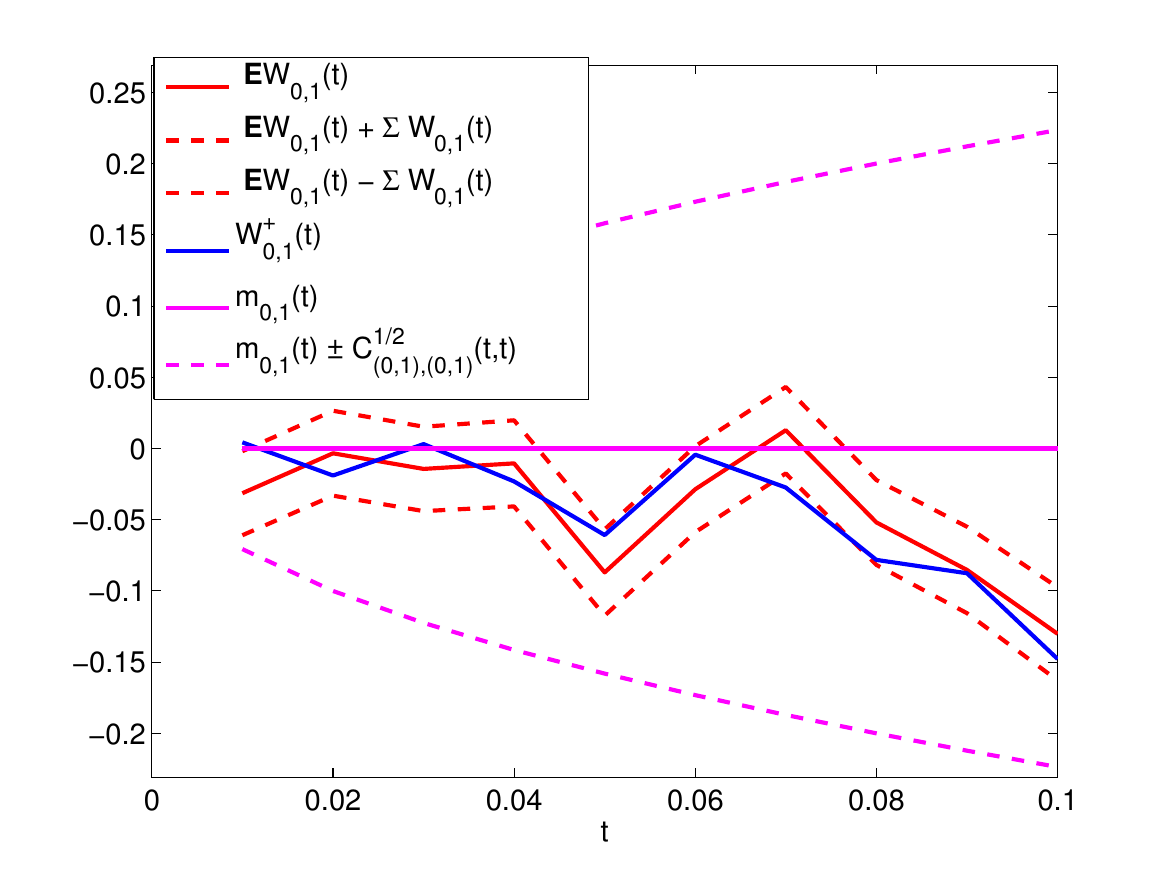}
 \includegraphics[width=0.5\textwidth]{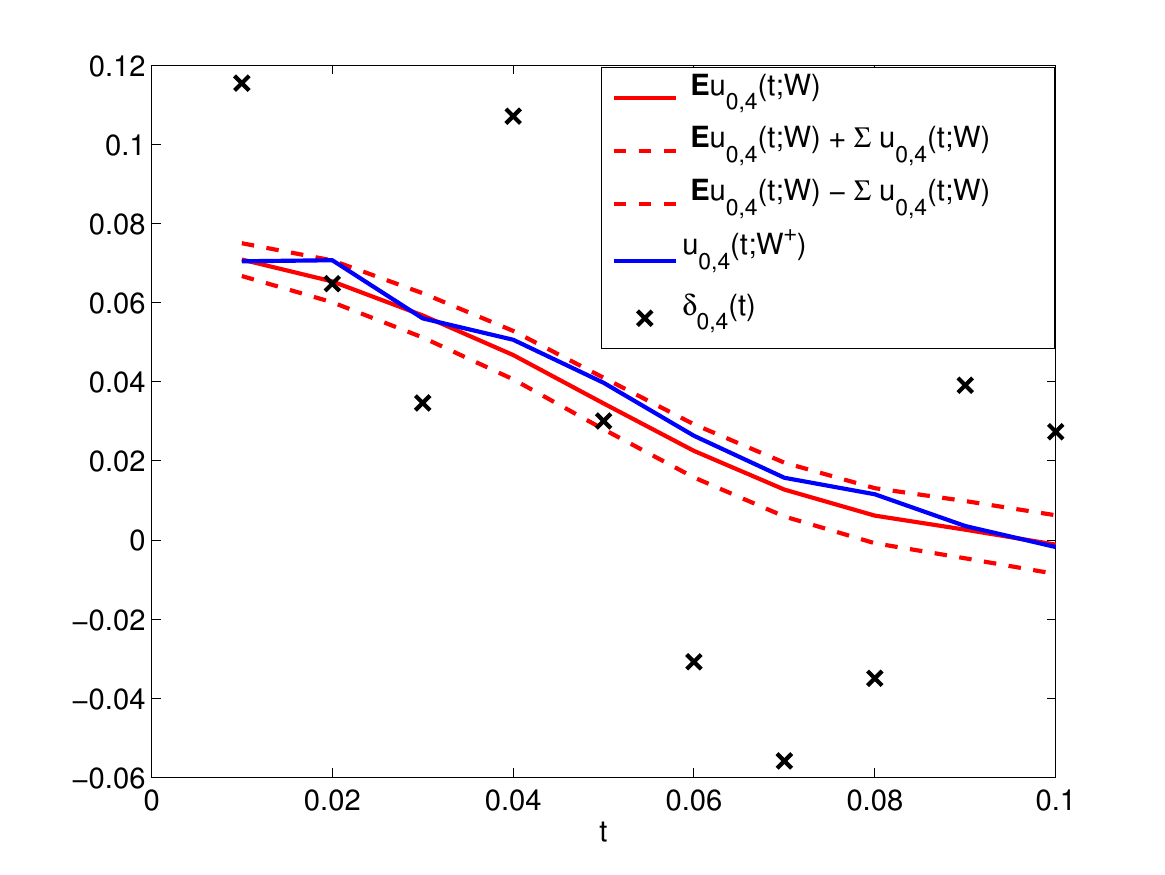}
 \includegraphics[width=0.5\textwidth]{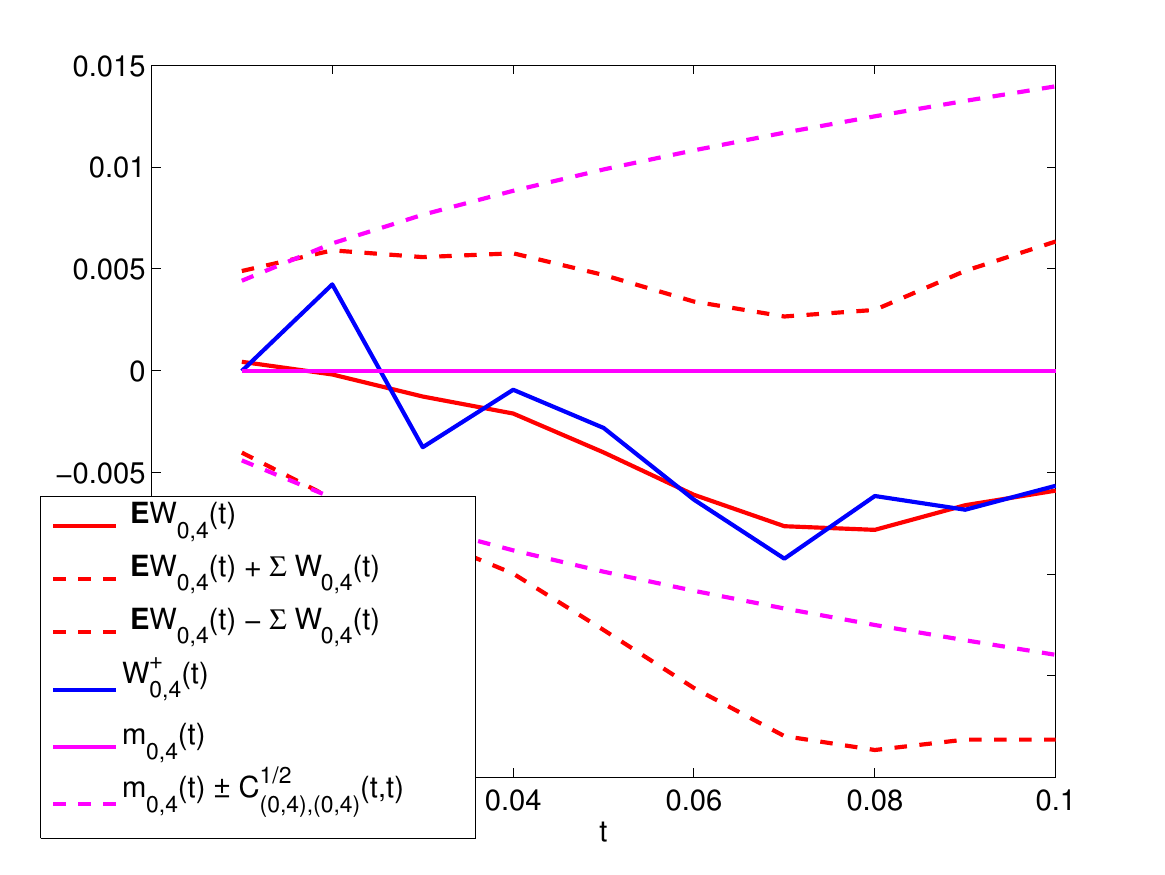} 
 \includegraphics[width=0.5\textwidth]{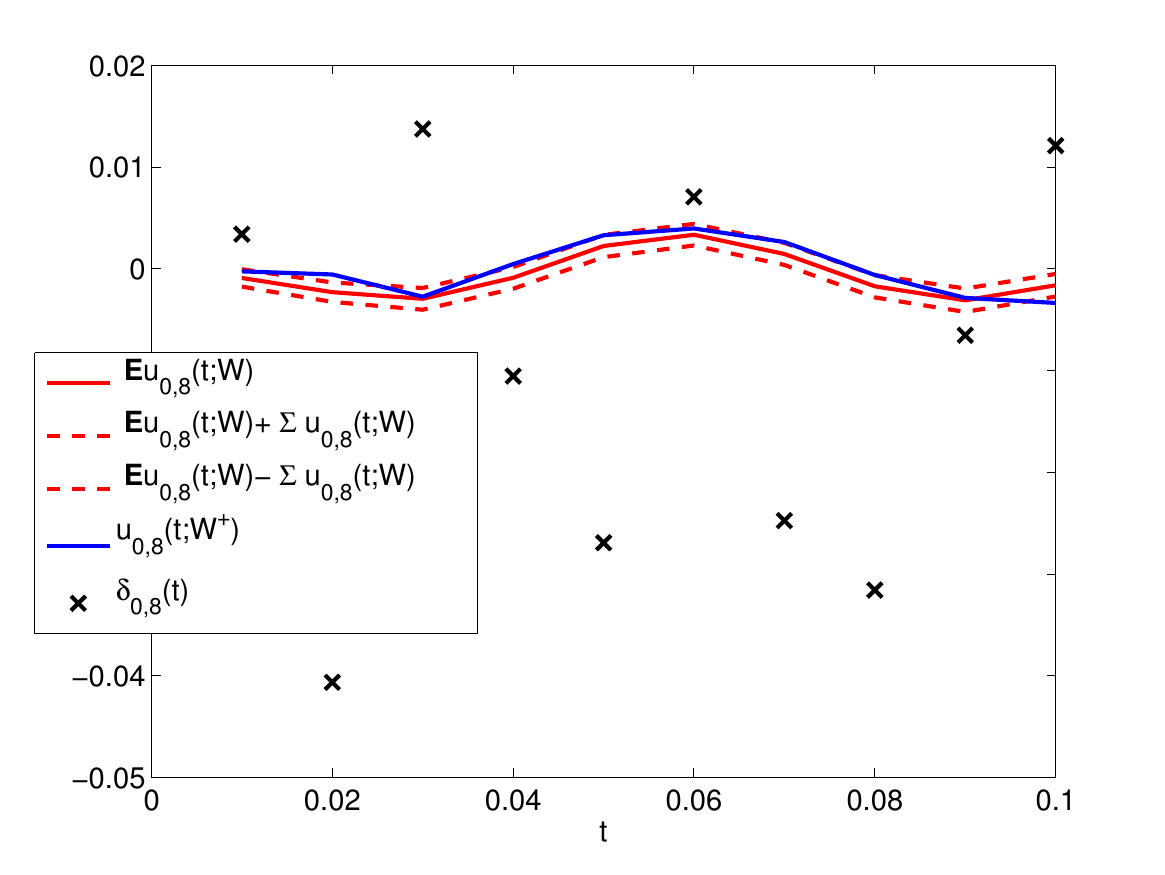}
 \includegraphics[width=0.5\textwidth]{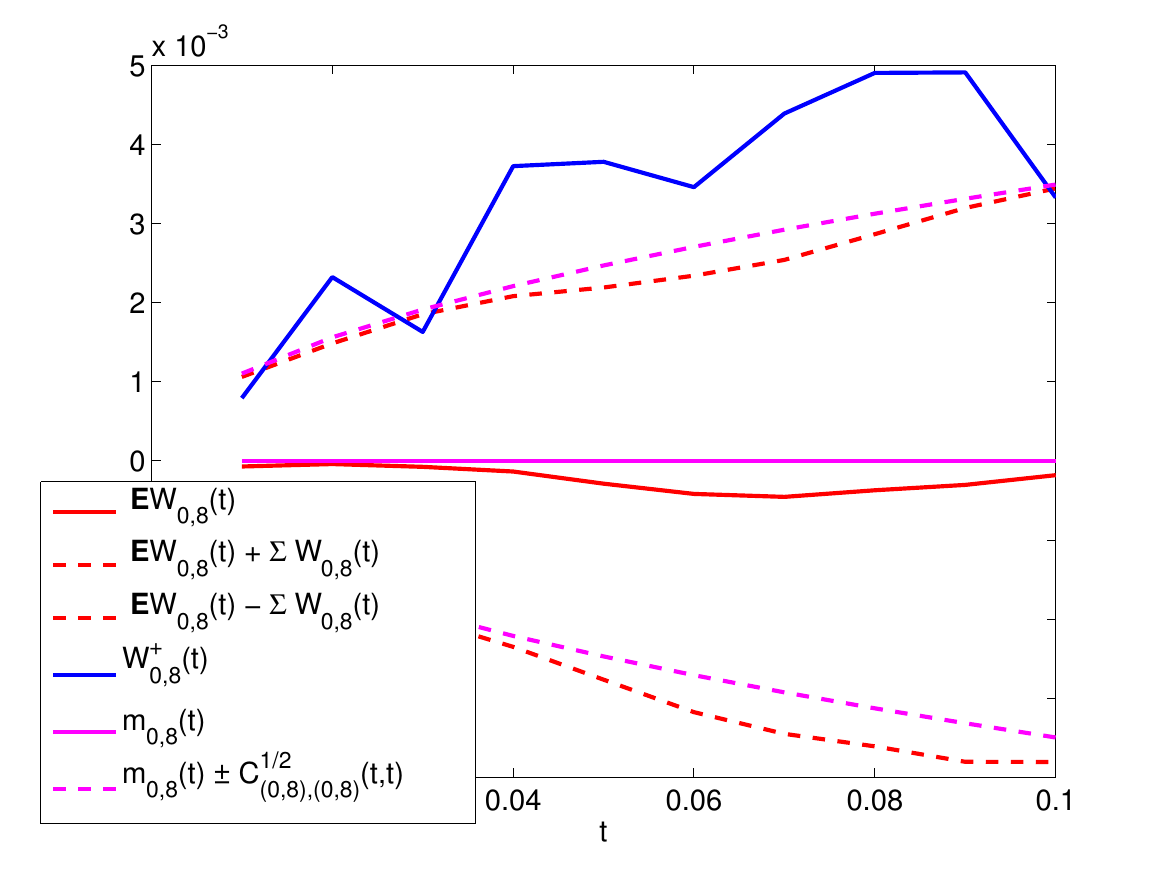}
 \caption{Full-field point wise observations. 
 The trajectories $u_{k}(t;W)$ (left) and 
$W_k$ (right), with $k=(0,1)$ (top),
$k=(0,4)$ (middle), and $k=(0,8)$ (bottom).  
Shown are expected values and 
standard deviation intervals as well as true values.
The right hand images also show the expected value 
and standard deviation of the prior, indicating the 
decreasing information content of the data for the 
increasing wave numbers.}
\label{traj}
\end{figure}

\begin{figure}
 \includegraphics[width=0.32\textwidth]{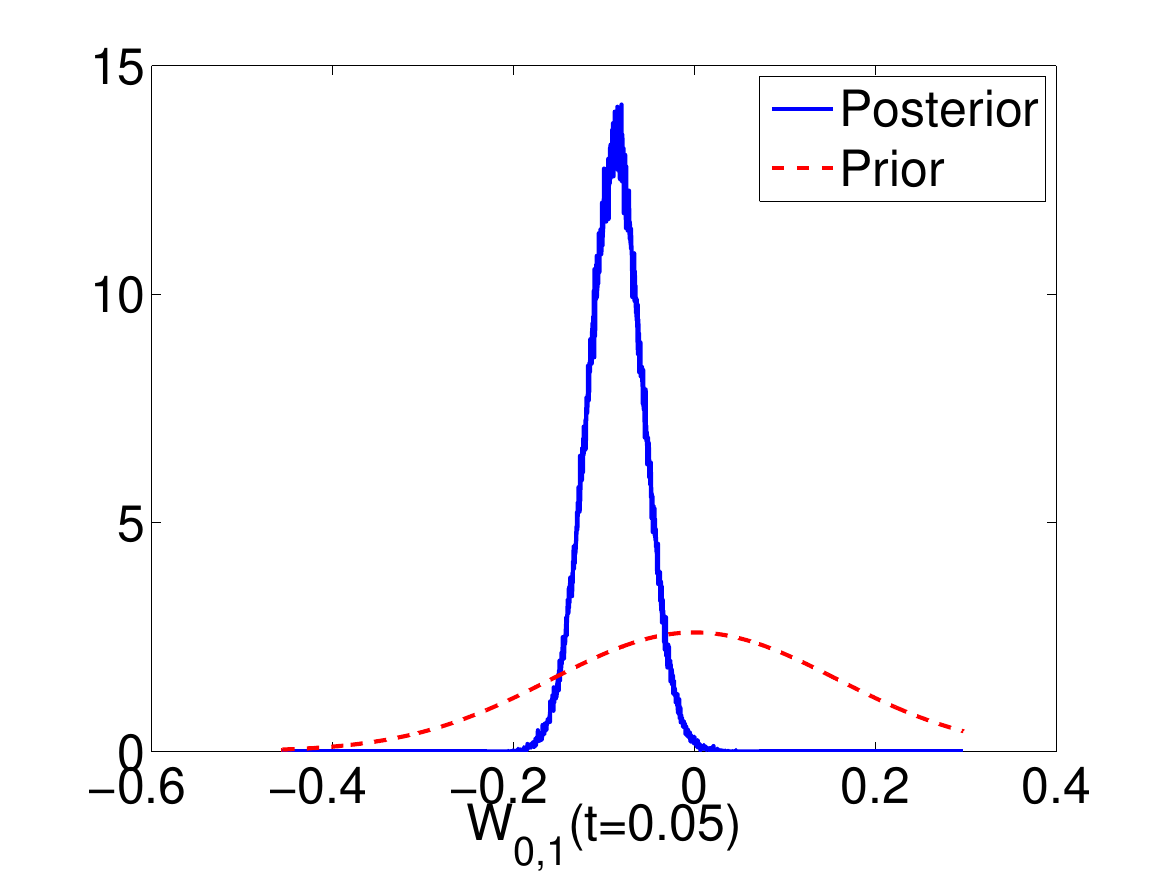}
 \includegraphics[width=0.32\textwidth]{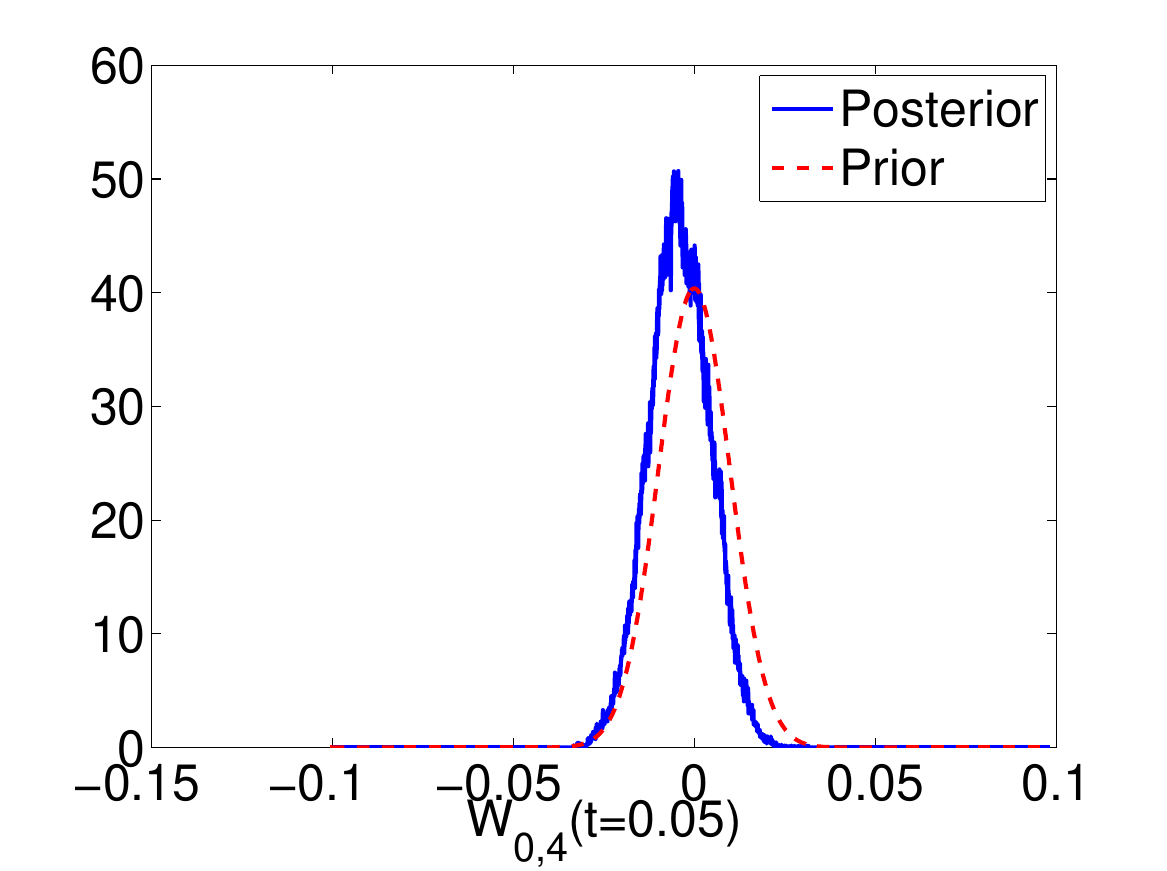}
 \includegraphics[width=0.32\textwidth]{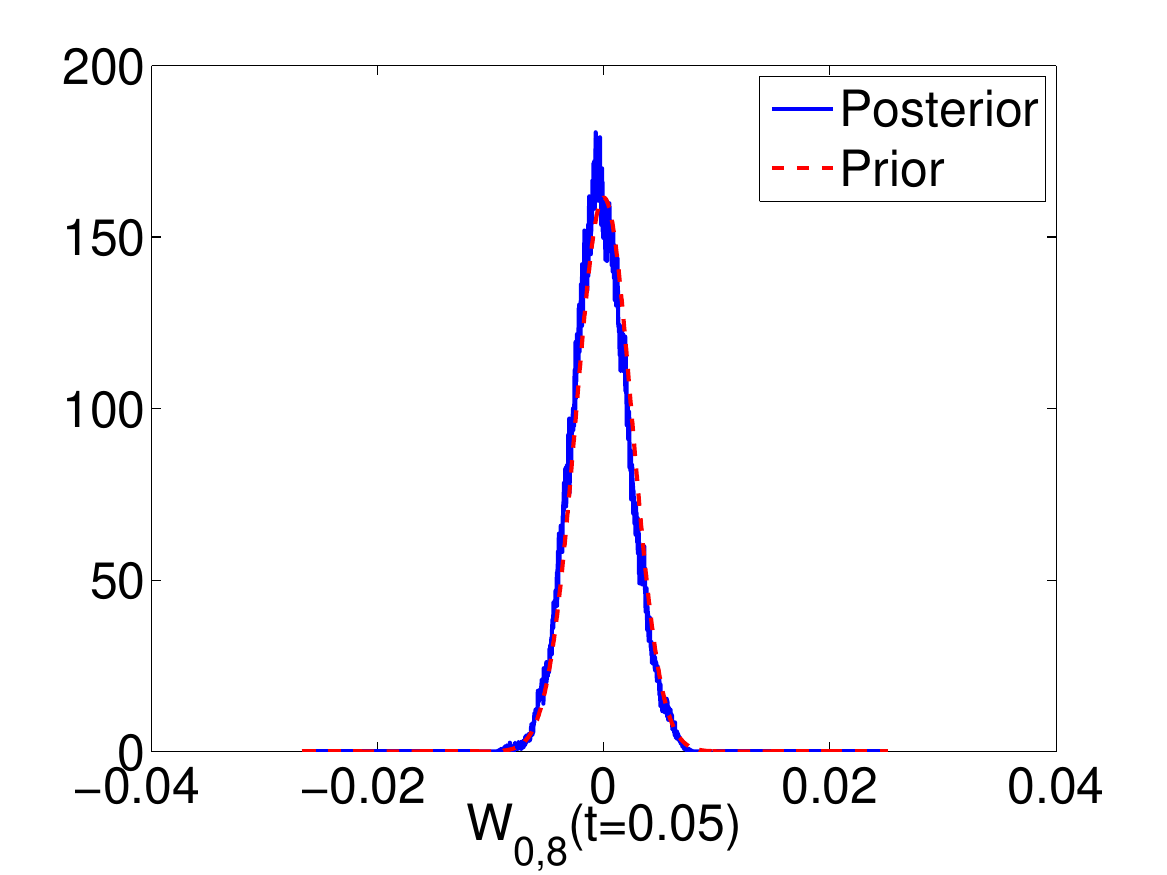}
 \caption{Full-field point wise observations. 
 The histograms of the posterior distribution of
 $W_k(t=0.05)$, for $k=(0,1)$ (left),
$k=(0,4)$ (middle), and $k=(0,8)$ (right).  
These plots again illustrate the 
decreasing information content of the data for the 
increasing wave numbers.}
\label{hist}
\end{figure}

We first consider the observation set-up (i) 
where pointwise observations of the entire velocity field are made.
The true initial and final 
conditions are plotted in Figure \ref{prof}, top two panels,
for the vorticity field $w$; the middle two panels of Figure \ref{prof}
show the posterior mean of the same quantities and indicate that
the data is fairly informative, since they closely resemble the truth;
the bottom two panels of Figure \ref{prof} show the absolute
difference between the fields in the top and middle panels.
The true trajectory, together with the posterior mean and 
one standard deviation interval around the mean, 
are plotted in Figure \ref{traj}, for the wavenumbers $(0,1)$, $(0,4)$,
and $(0,8)$, and for both the driving Brownian motion $W$ (right) and
the velocity field $u$ (left). This figure indicates that the data
is very informative about the $(0,1)$ mode, but less so concerning the
$(0,4)$ mode, and there is very little information in the $(0,8)$ mode.  
In particular for the $(0,8)$ mode the mean and standard deviation exhibit
behaviour similar to that under the prior whereas for the $(0,1)$ mode they
show considerable improvement over the prior in both position of the
mean and width of standard deviations.  {The posterior on the $(0,4)$ mode
has gleaned some information from the data as the mean has shifted
considerably from the prior; the variance remains similar to that
under the prior, however, so uncertainty in this mode has not been reduced.}
Figure \ref{hist}
shows the histograms of the prior and posterior for the same 3 modes
as in Fig. \ref{traj} at the center time $t=0.05$.  One can see here even
more clearly that the data is very informative about the $(0,1)$ mode 
in the left panel, less so but somewhat about the $(0,4)$ mode in the center panel, 
and it is not informative at all about the $(0,8)$ mode in the right panel.

\begin{figure}
 \includegraphics[width=0.5\textwidth]{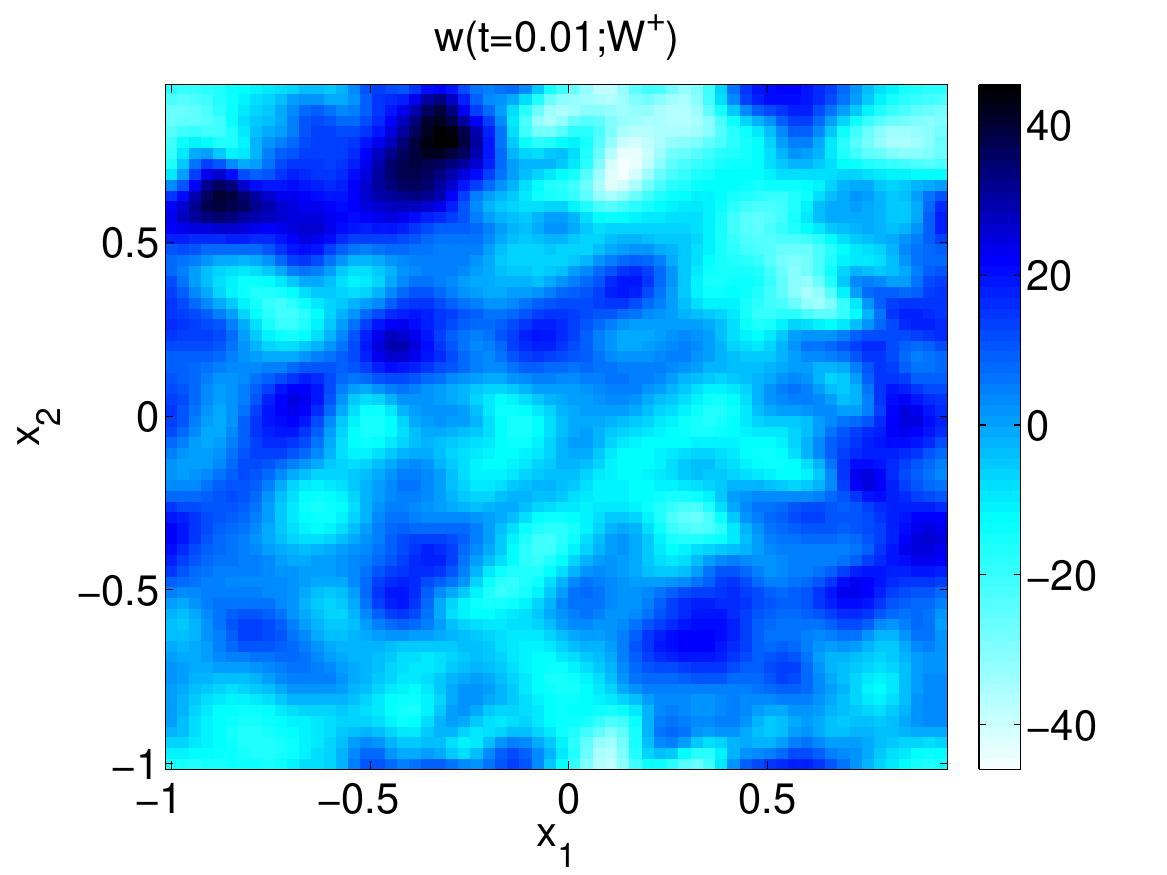}
 \includegraphics[width=0.5\textwidth]{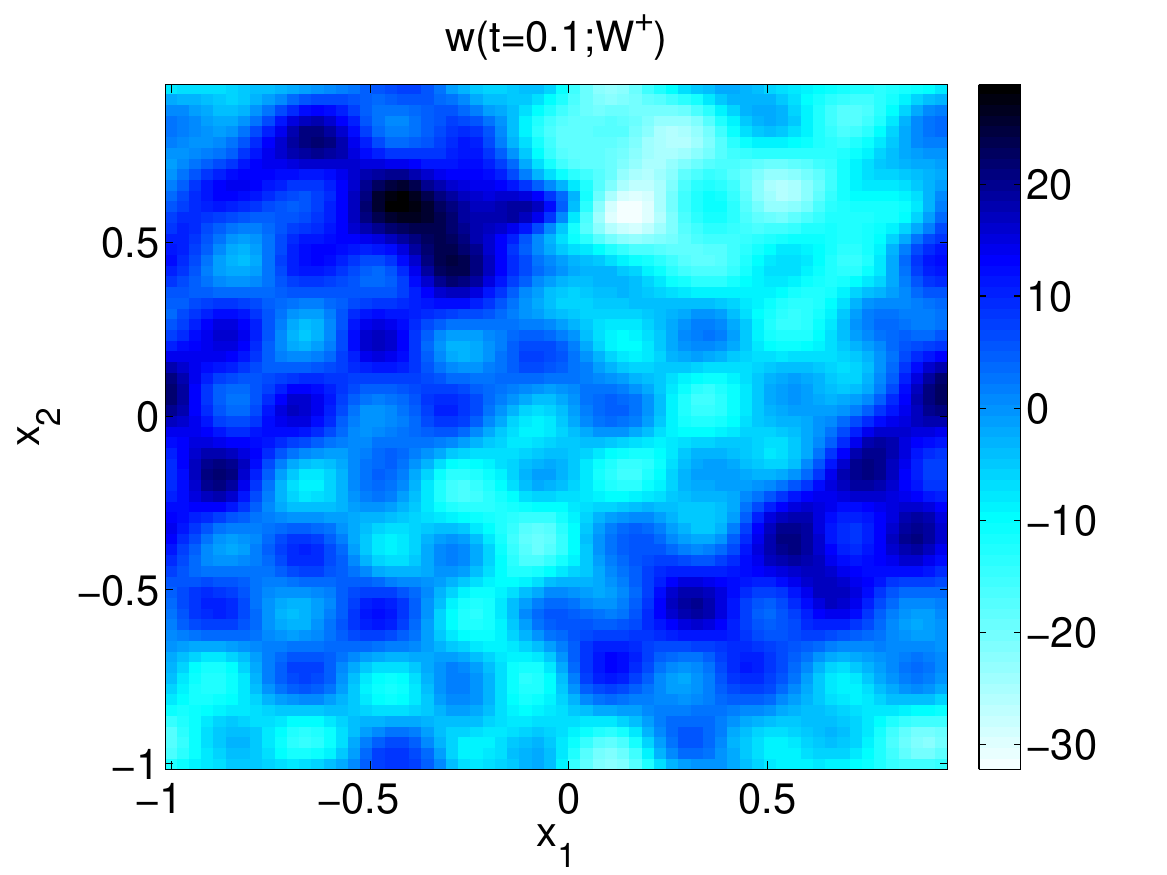}
 \includegraphics[width=0.5\textwidth]{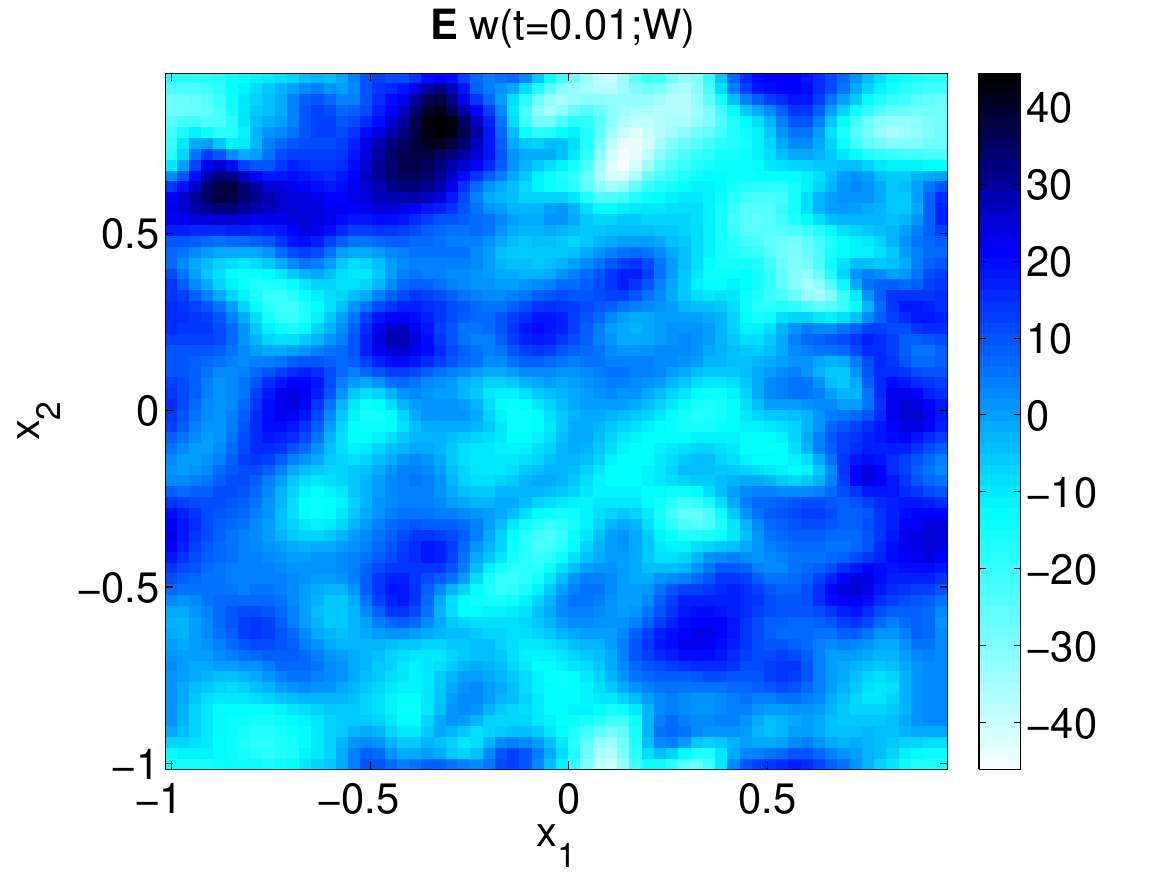}
 \includegraphics[width=0.5\textwidth]{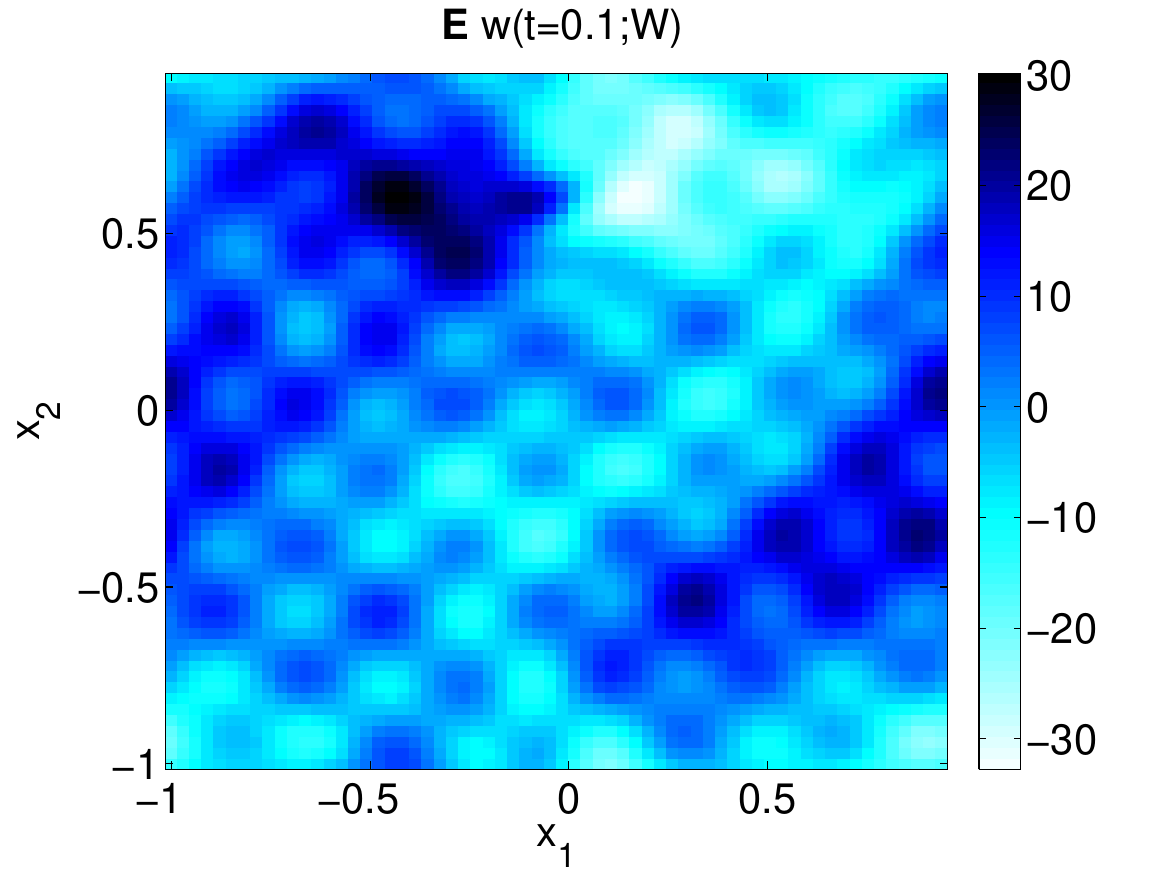}
 \includegraphics[width=0.5\textwidth]{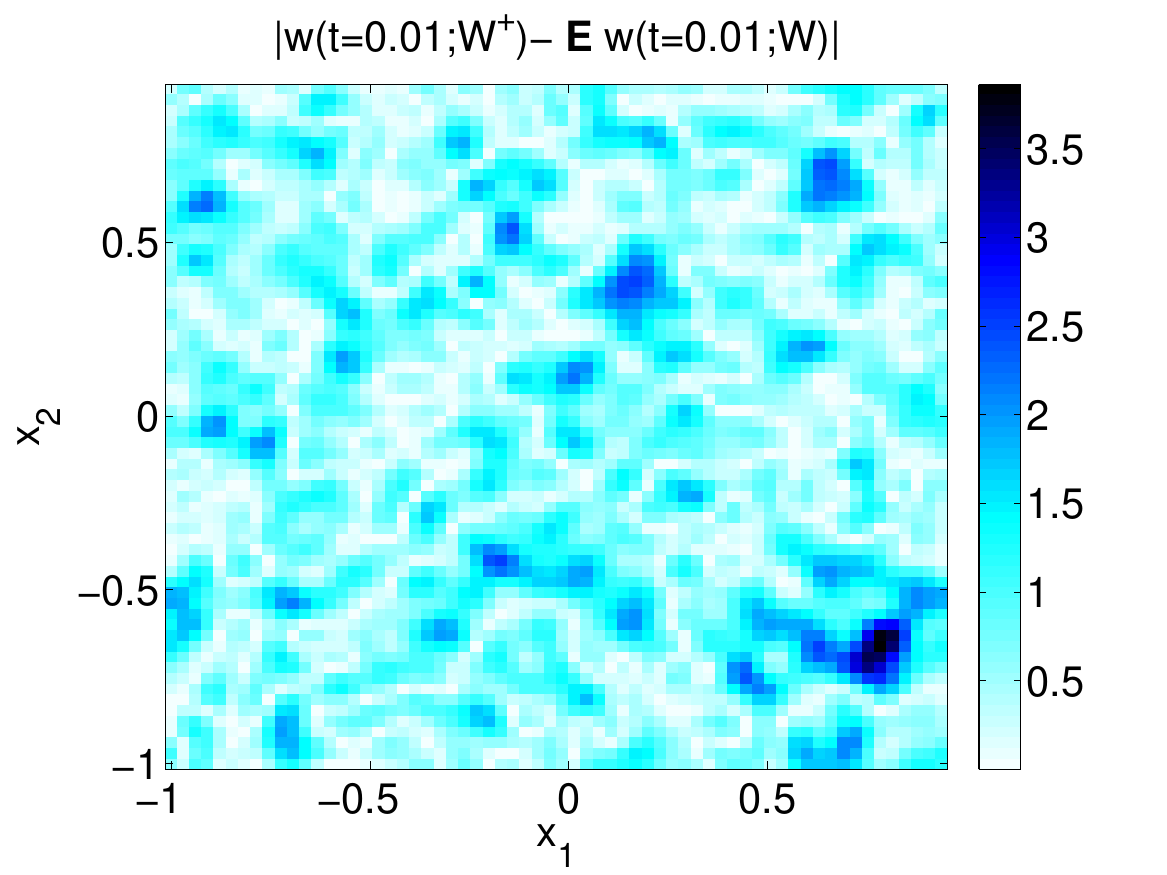}
 \includegraphics[width=0.5\textwidth]{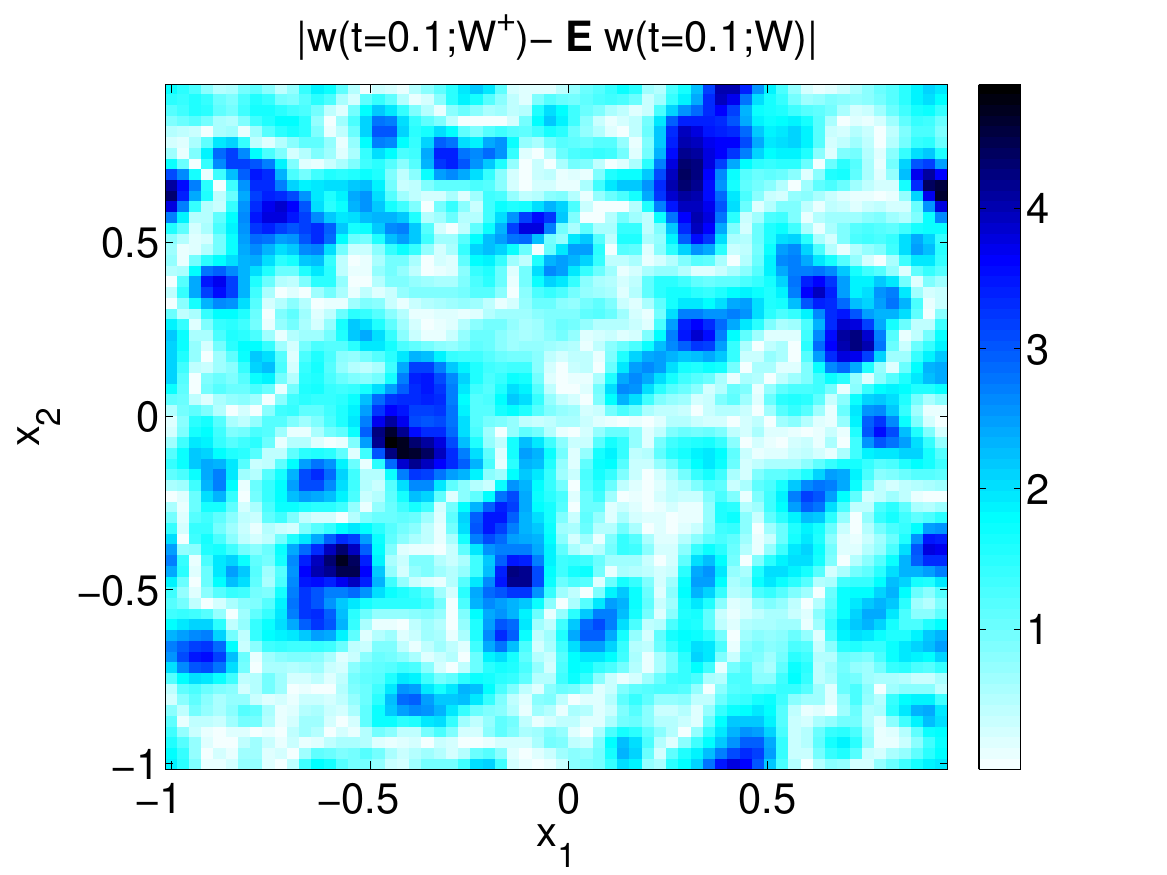} 
 \caption{Observation of Fourier modes $\{\phi_k\}_{|k|<4}$. 
 The truth (top), expected value (middle), and
 absolute distance between them (bottom) of the vorticity
 $w(t;W)$, for $t=0.01$ 
 (left, relative $L^2$ error $e=0.0044$) 
 and $t=0.1$ (right, $e=0.0249$). Notice the similarity to the results
 of Fig. \ref{prof}.}
\label{prof2}
\end{figure}

\begin{figure}
 \includegraphics[width=0.5\textwidth]{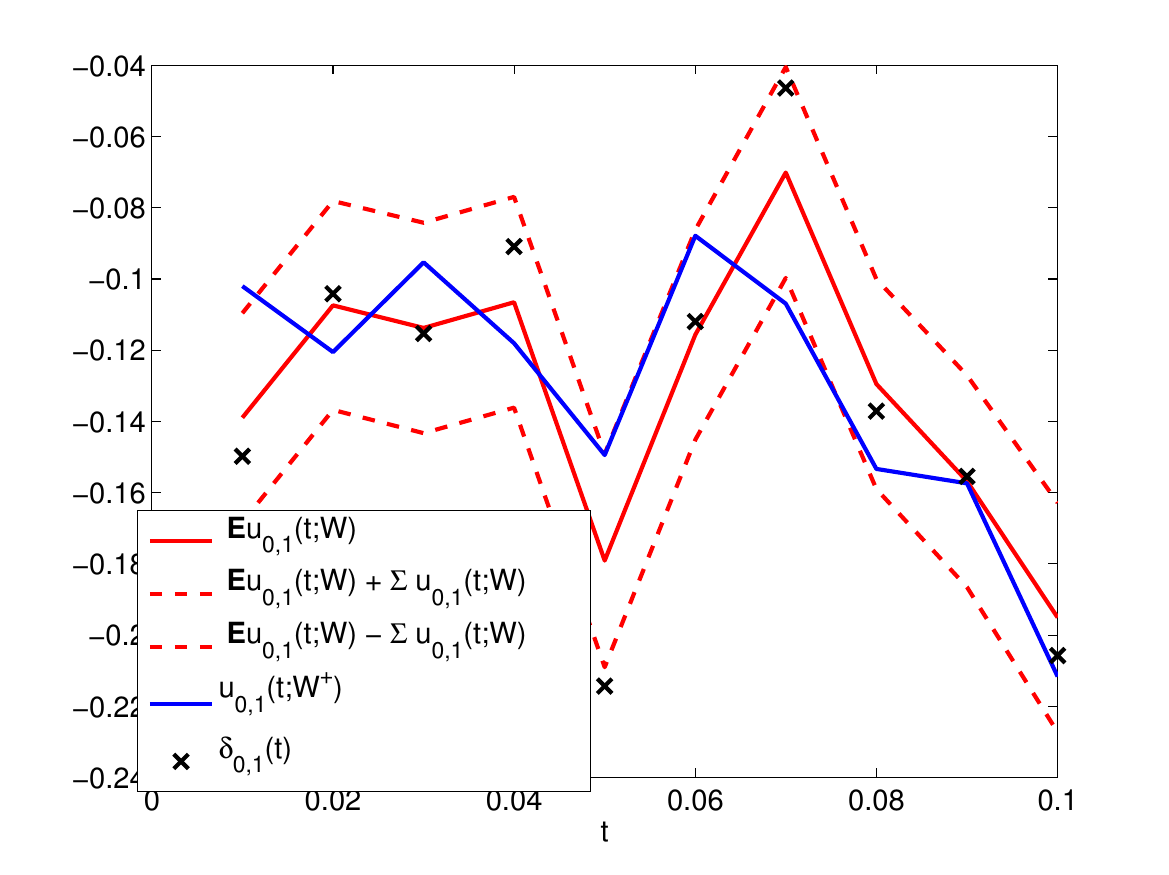}
 \includegraphics[width=0.5\textwidth]{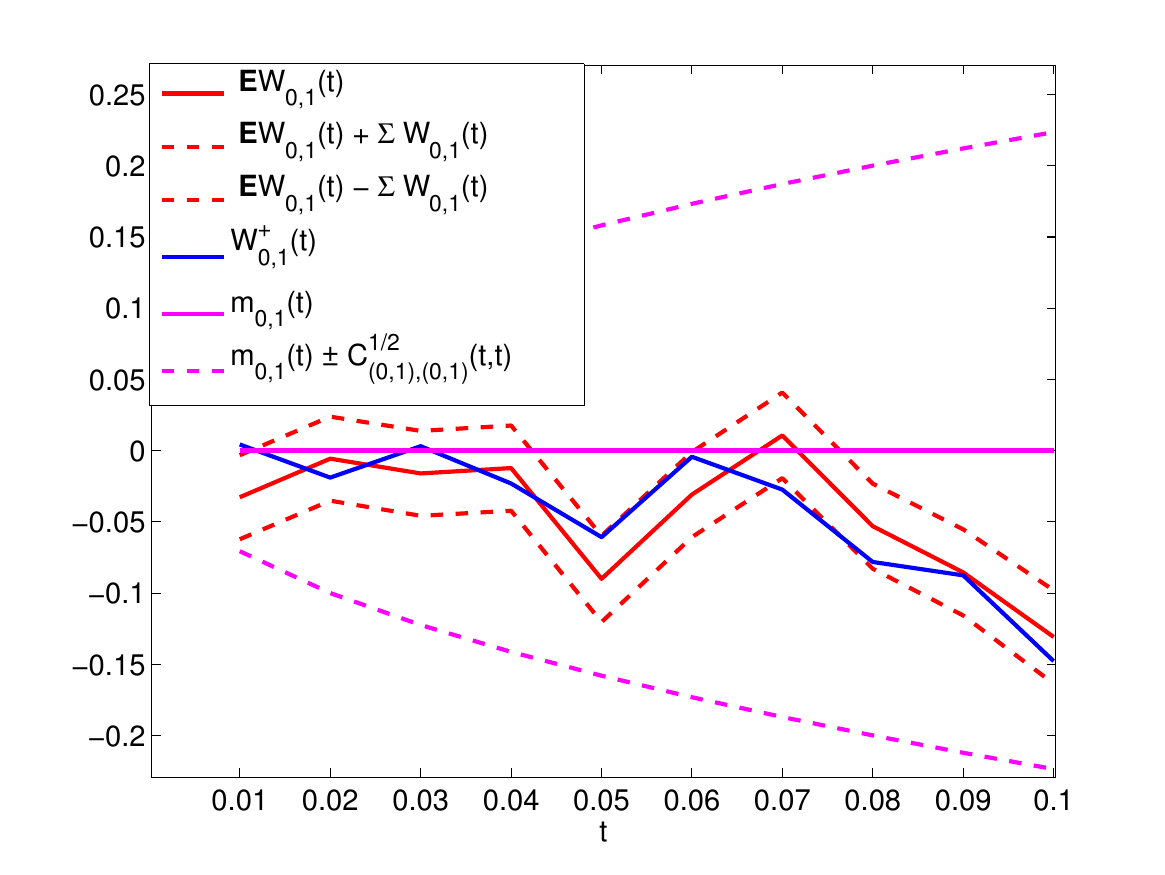}
 \includegraphics[width=0.5\textwidth]{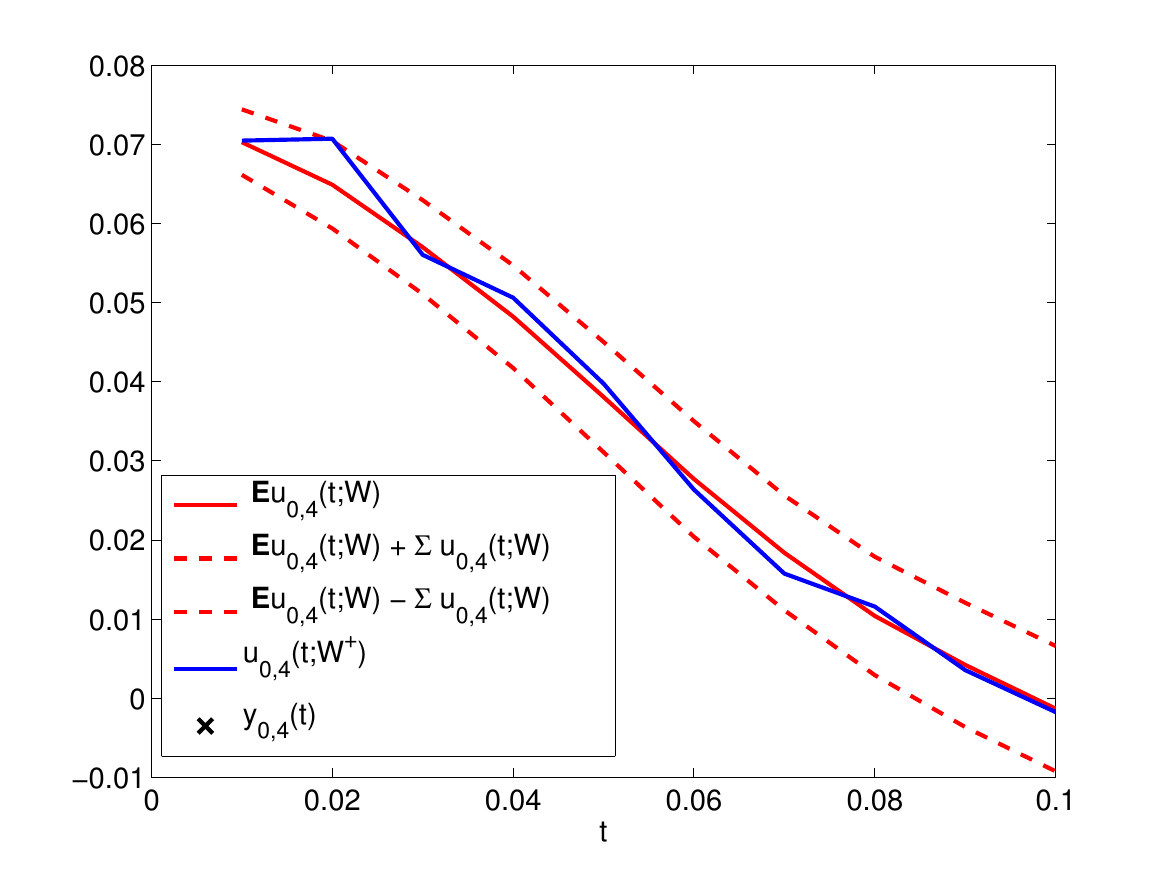}
 \includegraphics[width=0.5\textwidth]{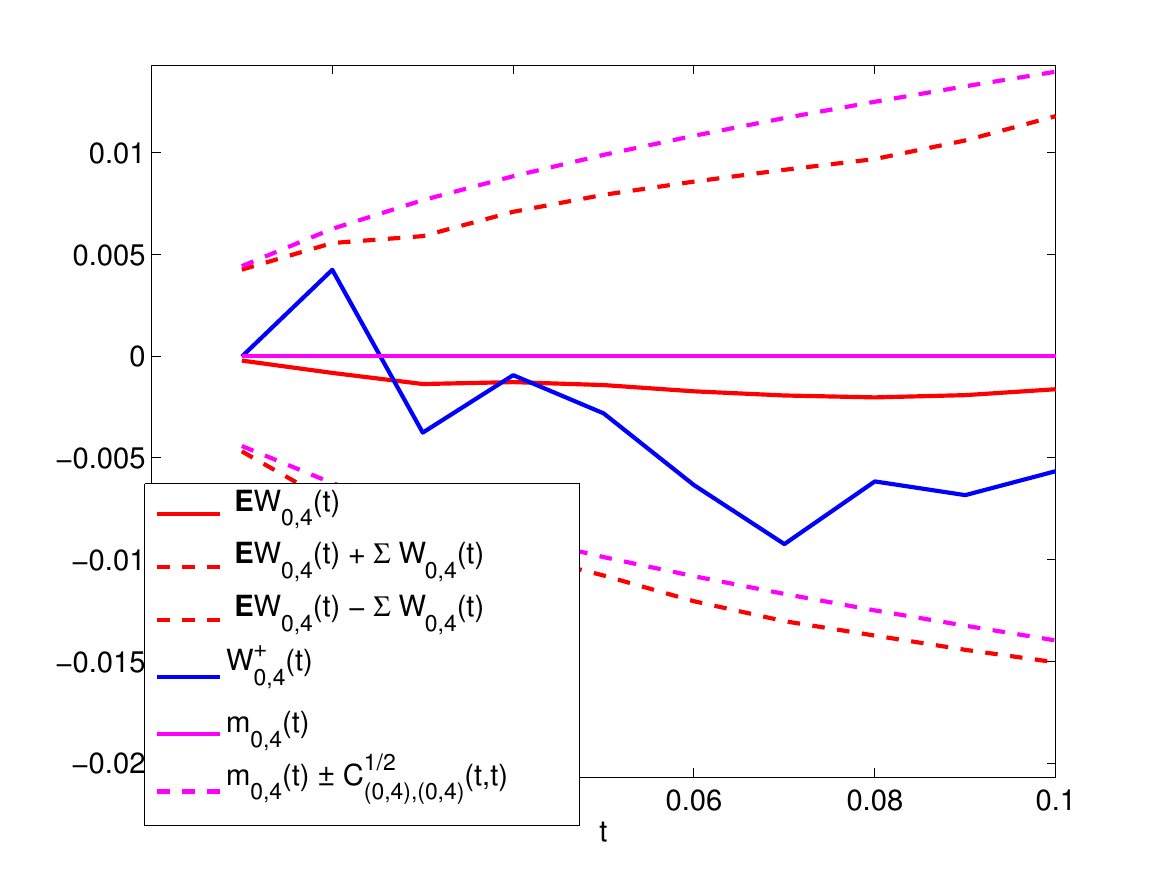} 
 \includegraphics[width=0.5\textwidth]{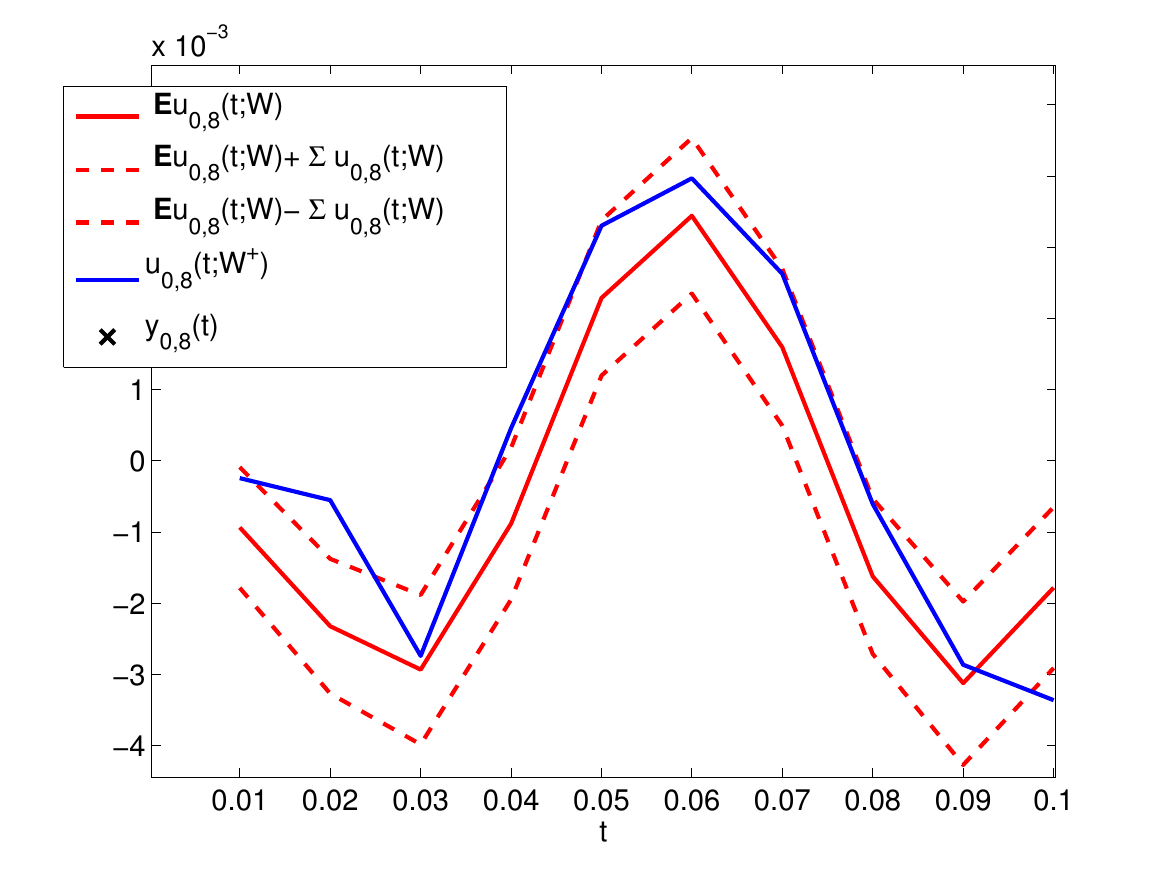}
 \includegraphics[width=0.5\textwidth]{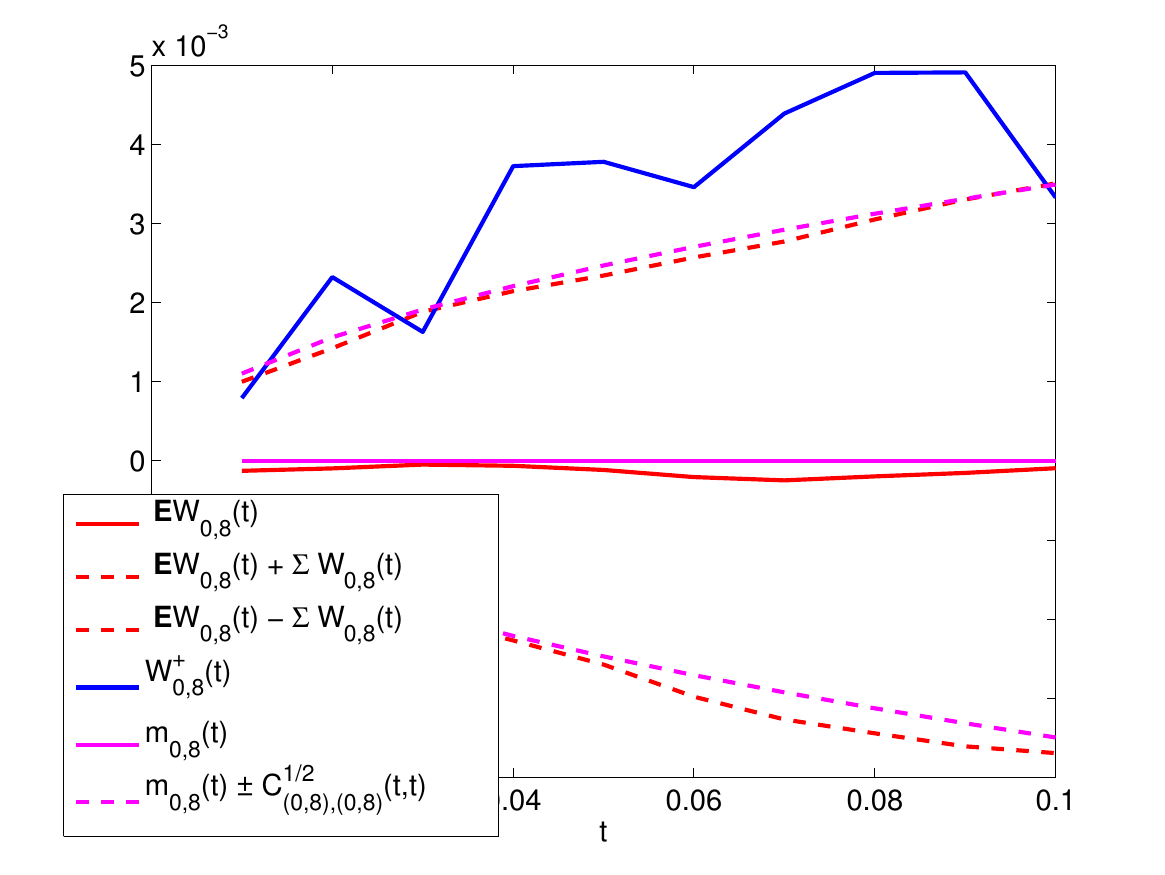}
 \caption{Observation of Fourier modes $\{\phi_k\}_{|k|<4}$.  
 The trajectories $u_{k}(t;W)$ (left) and 
$W_k$ (right), with $k=(0,1)$ (top),
$k=(0,4)$ (middle), and $k=(0,8)$ (bottom).  
Shown are expected values and 
standard deviation intervals as well as true values.
The right hand images also show the expected value 
and standard deviation of the prior, indicating the 
decreasing information content of the data for the 
increasing wave numbers.}
\label{traj2}
\end{figure}

\begin{figure}
 \includegraphics[width=0.32\textwidth]{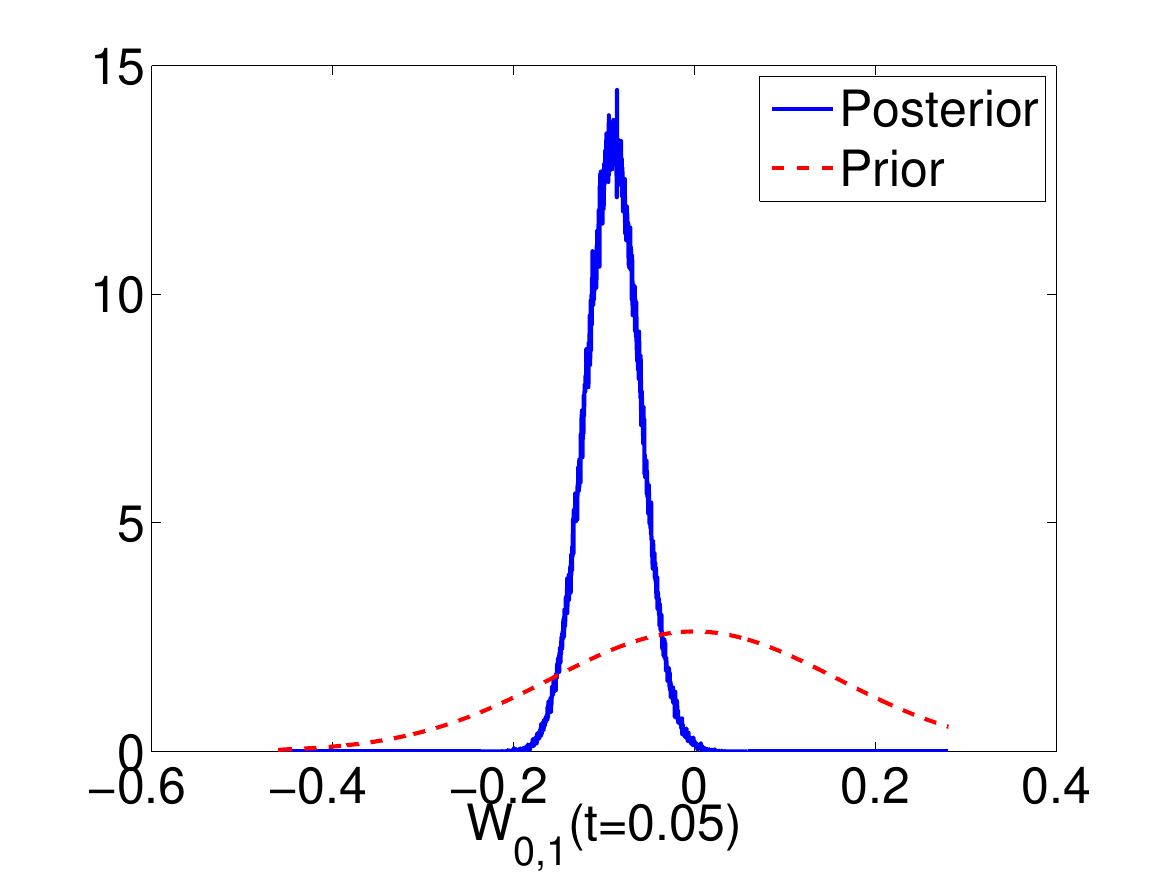}
 \includegraphics[width=0.32\textwidth]{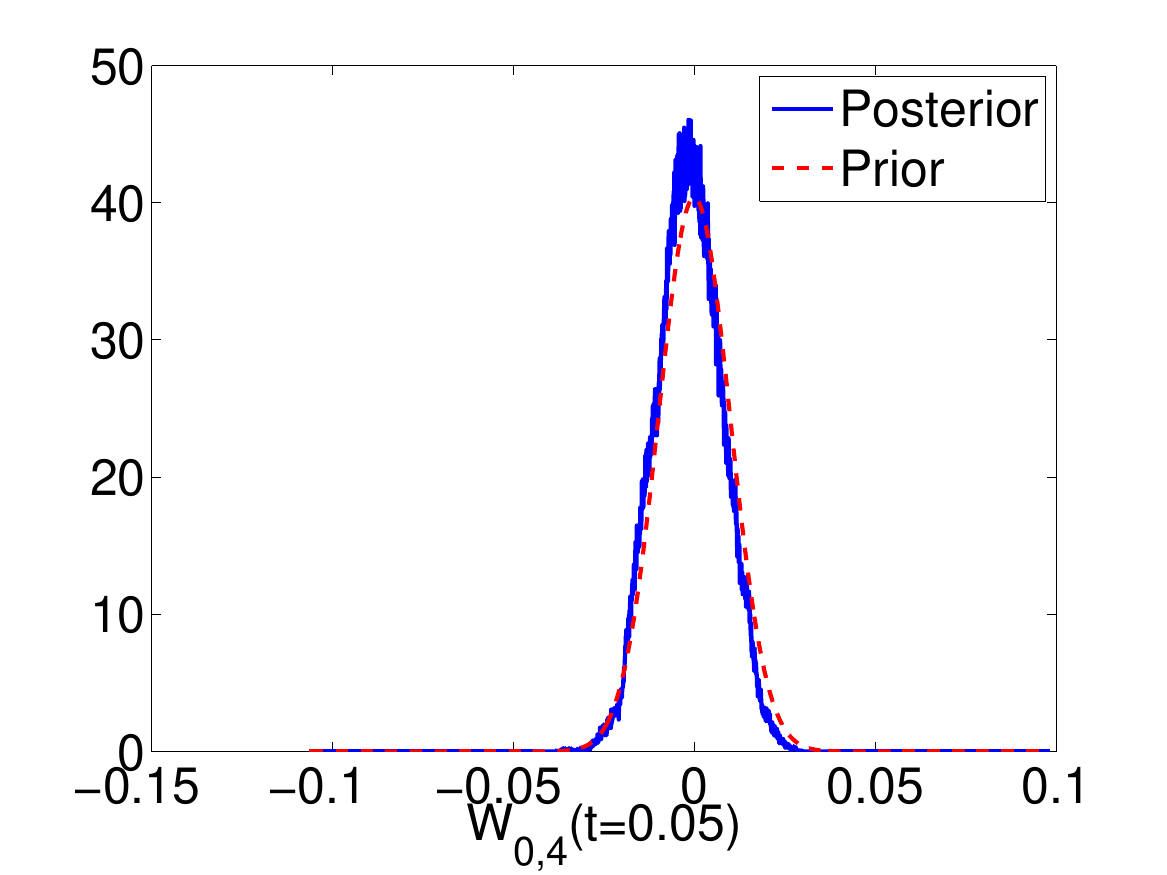}
 \includegraphics[width=0.32\textwidth]{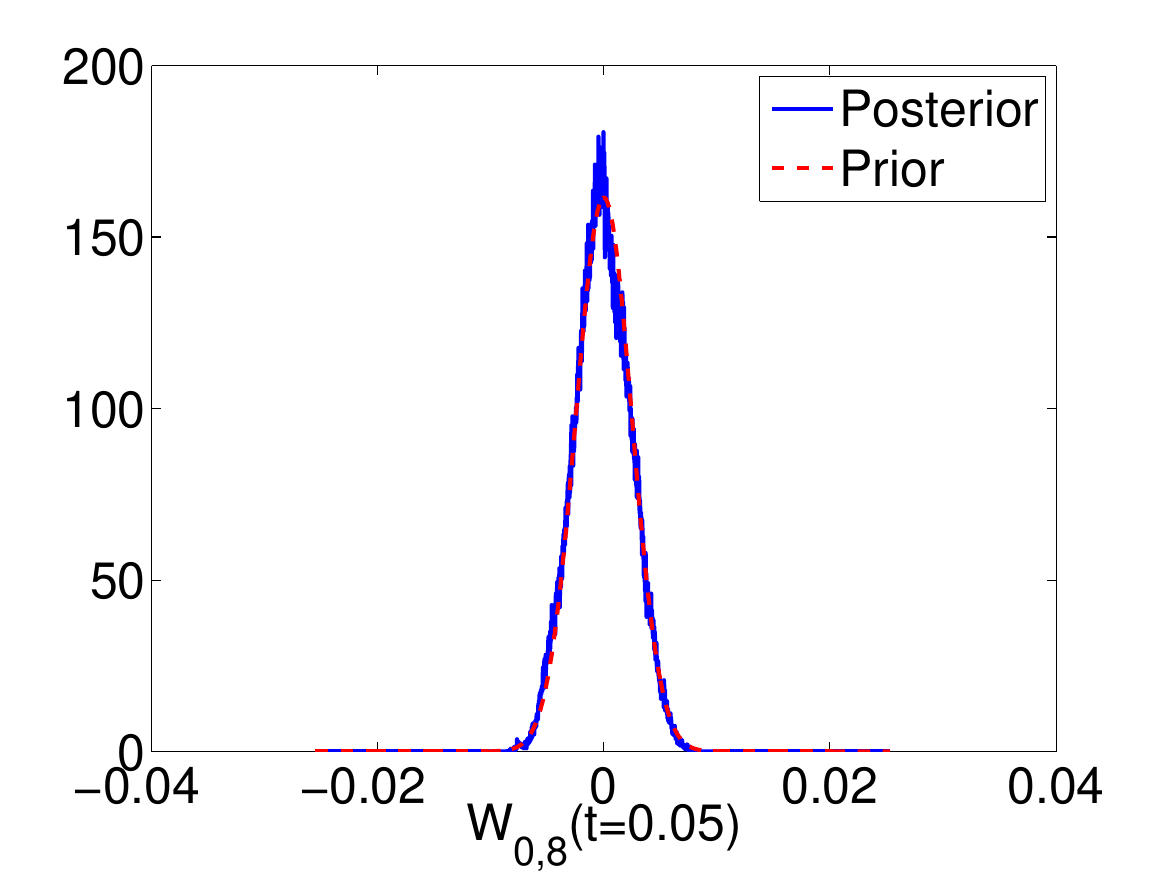}
 \caption{Observation of Fourier modes $\{\phi_k\}_{|k|<4}$.  
 The histograms of the posterior distribution of
 $W_k(t=0.05)$, for $k=(0,1)$ (left),
$k=(0,4)$ (middle), and $k=(0,8)$ (right).  
These plots again illustrate the 
decreasing information content of the data for the 
increasing wave numbers.  Notice the 
middle panel in which one notices the posterior on
$W_{0,4}(t=0.05)$ is much closer to the prior than
in Fig. \ref{hist}.}
\label{hist2}
\end{figure}

Figures  \ref{prof2}, \ref{traj2}, and \ref{hist2} are the same as 
Figures \ref{prof}, \ref{traj}, and \ref{hist} except for the case of (ii) observation
of low Fourier modes.  Notice that the difference in the spatial fields are
difficult to distinguish by eye, and indeed the relative errors even 
agree to threshold $10^{-3}$.  However, we can see that now the unobserved
$(0,4)$ mode in the center panels of Figs. \ref{traj2} and \ref{hist2} is not 
informed by the data and remains distributed approximately like the prior.

\begin{acknowledgement} VHH is supported by a start-up grant from Nanyang Technological University, 
AMS is grateful to EPSRC, ERC, ESA and ONR for financial support for this work, 
and KJHL is grateful to ESA and the SRI-UQ Center at KAUST for financial support.
\end{acknowledgement}
\bibliographystyle{plain}
\bibliography{NSE}

\end{document}